\documentclass[preprint,5p, twocolumn]{elsarticle}

\usepackage{amsthm}
\usepackage{amsmath}
\usepackage{amssymb}
\usepackage{dsfont}
\usepackage{tikz}
\usetikzlibrary{patterns,arrows, math}
\usetikzlibrary{shapes, arrows}
\tikzstyle{block} = [draw, fill=white, rectangle, minimum height=3em, minimum width=6em]
\tikzstyle{output} = [coordinate]
\tikzstyle{input} = [coordinate]
\usepackage{pgfplots}
\pgfplotsset{compat=1.18} 

\usepackage{textcomp}   
\def\BibTeX{{\rm B\kern-.05em{\sc i\kern-.025em b}\kern-.08em
    T\kern-.1667em\lower.7ex\hbox{E}\kern-.125emX}}

\usepackage[colorlinks=true]{hyperref}
\usepackage{float}
\usepackage{subcaption}
\usepackage{comment}
\usepackage{enumerate}

\usepackage[normalem]{ulem}
\usepackage{mathtools}

\usepackage{appendix}
\usepackage{cleveref}

\crefname{appsec}{Appendix}{Appendices}

\usepackage{algorithm}
\usepackage[italicComments=true, beginComment=//~]{algpseudocodex}

\usepackage{nicefrac}

\newtheorem{theorem}{Theorem}[section]
\newtheorem{assumption}{Assumption}
\newtheorem{definition}{Definition}[section]
\newtheorem{remark}{Remark}[section]
\newtheorem{lemma}{Lemma}[section]

\Crefname{definition}{Definition}{Definitions}
\Crefname{assumption}{Assumption}{Assumptions} 
\Crefname{lemma}{Lemma}{Lemmata}
\Crefname{corollary}{Corollary}{Corollaries}
\Crefname{remark}{Remark}{Remarks}
\Crefname{theorem}{Theorem}{Theorems}

\makeatletter
\newcommand\RedeclareMathOperator{%
  \@ifstar{\def\rmo@s{m}\rmo@redeclare}{\def\rmo@s{o}\rmo@redeclare}%
}
\newcommand\rmo@redeclare[2]{%
  \begingroup \escapechar\m@ne\xdef\@gtempa{{\string#1}}\endgroup
  \expandafter\@ifundefined\@gtempa
     {\@latex@error{\noexpand#1undefined}\@ehc}%
     \relax
  \expandafter\rmo@declmathop\rmo@s{#1}{#2}}
\newcommand\rmo@declmathop[3]{%
  \DeclareRobustCommand{#2}{\qopname\newmcodes@#1{#3}}%
}
\@onlypreamble\RedeclareMathOperator
\makeatother

{\left(\begin{smallmatrix}}%
{\end{smallmatrix}\right)}%

{\left[\begin{smallmatrix}}%
{\end{smallmatrix}\right]}%

\newcommand{\N}{\mathds{N}}
\newcommand{\R}{\mathds{R}}
\newcommand{\Rp}{\R_{\geq0}}

\newcommand{\Impl}{\Longrightarrow}

\newcommand{\fa}{\ \forall \, }
\newcommand{\ex}{\ \exists \, }

\newcommand{\al}{\left \langle}
\newcommand{\ar}{\right \rangle}

\newcommand{\nl}{\left\|}
\newcommand{\nr}{\right\|}
\newcommand{\cbl}{\left\lbrace }
\newcommand{\cbr}{\right\rbrace }

\newcommand{\Norm}[2][ ]{\nl #2 \nr_{#1}}
\newcommand{\SNorm}[1]{\Norm[\infty]{#1}}

\newcommand{\setdef}[2]{\cbl\ #1\ \left|\ \vphantom{#1} #2\ \right.\cbr}

\newcommand{\cB}{\mathcal{B}}
\newcommand{\cC}{\mathcal{C}}
\newcommand{\cD}{\mathcal{D}}

\newcommand{\cU}{\mathcal{U}}

\newcommand{\cT}{\mathcal{T}}

\newcommand{\cF}{\mathcal{F}}
\newcommand{\cG}{\mathcal{G}}
\newcommand{\cN}{\mathcal{N}}
\newcommand{\cY}{\mathcal{Y}}

\newcommand{\cX}{\mathcal{X}}

\newcommand{\oT}{\mathbf{T}}

\newcommand{\Lip}{\rm{Lip}}

\DeclareMathOperator*{\rf}{ref}
\DeclareMathOperator*{\esssup}{ess\,sup}

\DeclareMathOperator*{\loc}{loc}

\newcommand{\ve}{\varepsilon}

\renewcommand{\l}{\lambda}

\newcommand{\con}{\mathcal{C}}

\RedeclareMathOperator*{\Im}{Im}
\RedeclareMathOperator*{\Re}{Re}
\renewcommand{\phi}{\varphi}

\renewcommand{\d}{\ \text{d}}

\newcommand{\vp}{\varphi}

\newcommand{\dd}[2][ ]{\tfrac{\text{\normalfont d}#1}{\text{\normalfont d}#2}}

\newcommand{\fM}{f_{\rm max}}
\newcommand{\gM}{g_{\rm max}}
\newcommand{\gm}{g_{\rm min}}

\journal{Systems \& Control Letters}

\begin{document}

\ifpdf
\hypersetup{
  pdftitle={Sampled-data funnel control and its use for safe continual learning}
}
\fi

\begin{frontmatter}
\title{
    Sampled-data funnel control and its use for safe continual learning\tnoteref{Support}\tnotetext[Support]{
Funding: L.~Lanza, D.~Dennstädt, and K.~Worthmann gratefully acknowledge funding by the Deutsche Forschungsgemeinschaft (DFG, German Research Foundation; Project-IDs 471539468 and 507037103).
L.~Lanza and P.~Schmitz are 
grateful for 
support by the Carl Zeiss Foundation (VerneDCt -- Project No.\ 2011640173 and DeepTurb--- Project No.\ 2018-02-001).
G.\ D.\ \c{S}en gratefully acknowledges funding by the German Academic Exchange Service (DAAD).
}
}

\author[label_tuilm]{Lukas Lanza}
\ead{lukas.lanza@tu-ilmenau.de}

\author[label_tuilm,label_upb]{Dario Dennstädt}
\ead{dario.dennstaedt@uni-paderborn.de}

\author[label_tuilm]{Karl Worthmann}
\ead{karl.worthmann@tu-ilmenau.de}

\author[label_tuilm]{Philipp Schmitz}
\ead{philipp.schmitz@tu-ilmenau.de}

\author[label_itu,label_btu]{Gökçen Devlet \c{S}en\texorpdfstring{\corref{cor}}{}}
\ead{sen@b-tu.de}

\author[label_rug]{Stephan Trenn}
\ead{s.trenn@rug.nl}

\author[label_tuilm]{Manuel Schaller}
\ead{manuel.schaller@tu-ilmenau.de}

\affiliation[label_tuilm]{organization={Technische Universität Ilmenau, Institute of Mathematics, Optimization-based Control Group},%
            city={Ilmenau},
            country={Germany}}

\affiliation[label_upb]{organization={Universität Paderborn, Institut für Mathematik},%
            city={Paderborn},
            country={Germany}}

\affiliation[label_rug]{organization={University of Groningen, Bernoulli Institute for Mathematics, Computer
Science, and Artificial Intelligence},%
            city={Groningen},
            country={The Netherlands}}
\affiliation[label_itu]{organization={Istanbul Technical University, Department of Control and Automation Engineering},%
            city={Istanbul},
            country={Turkey}}
\affiliation[label_btu]{organization={Brandenburgische Technische Universität Cottbus-Senftenberg, Fachgebiet Regelungssysteme und Netzleittechnik},%
            city={Cottbus},
            country={Germany}}

\cortext[cor]{Corresponding Author: Gökçen Devlet \c{S}en}

\begin{abstract}
We propose a novel sampled-data output-feedback controller for nonlinear 
systems of arbitrary relative degree 
that ensures reference tracking within prescribed
error bounds.
We provide explicit bounds on the maximum input signal and the required uniform sampling time. A key strength of this approach is its capability to serve as a safety filter for various learning-based controller designs, enabling the use of learning techniques in safety-critical applications. We 
illustrate its versatility by integrating it with two different controllers: a reinforcement learning controller and a non-parametric predictive controller based on Willems et al.'s fundamental lemma.
Numerical simulations illustrate effectiveness of the combined controller design.
\end{abstract}


\end{frontmatter}

\section{Introduction}
In the context of output-reference tracking, \emph{funnel control} is an established adaptive high-gain control methodology, which guarantees satisfaction of a-priori fixed, possibly time-varying output constraints. Apart from imposing structural assumptions such as known relative degree, a high-gain property, and a bounded-input-bounded-state property of the internal dynamics, no system knowledge is required, see~\cite{BergIlch21} 
and the references therein. 
Pivotal for its functioning is the availability of the system's output as a time-continuous signal and the ability to continuously adapt the input signal. This requirement, however, is challenged by 
digital measurement devices and controllers. \\
Although funnel control has been successfully implemented in a sampled-data system with Zero-order
Hold (ZoH) for a sufficiently small sampling time~in~\cite{BergKaes20}, we are not aware of any
results rigorously showing that the output signal stays within the prescribed boundaries for ZoH funnel
control. 
In this paper, we address this disparity by proposing
a novel sampled-data feedback controller with ZoH.
The proposed controller ensures output tracking of a given reference signal within prescribed,  
possibly time-varying performance bounds -- at every time instant meaning that also the intersampling behavior is fully taken into account.
Balancing the need for a sufficiently large feedback gain for output tracking
and avoidance of overshooting (which
could violate error bounds
within one sampling period),
we derive uniform bounds on sampling rates and control inputs 
such that the imposed output constraints are satisfied along the closed loop leveraging coarse bounds
on the
system dynamics.
To the best of our knowledge, in funnel control uniform bounds on the input signal are only known if
the region of feasible initial values is further restricted \textit{and} the dynamics are known~\cite{BergIlch21}. 
While there have been several attempts to deal with the closely
related issue of input
saturation~\cite{Berg22,HuTren22,IlchTren04} and bang-bang
controller designs~\cite{LibeTren13b} exhibiting similarities to our approach, an analysis of combining a
ZoH with funnel control has not been conducted.

The controller proposed in this article includes an ``activation threshold'' to set the input to zero for small
tracking errors, akin to approaches in~\cite{Schenato09} and in~\cite{BergDenn23} using an
activation function, the $\lambda$-tracker~\cite{IlchRyan94}, or more broadly event- and self-triggered
controller designs, see e.g.~\cite{Heemels2021} and references therein.
This opens up the possibility for the controller to act as a safety filter for 
data-driven approaches
and (online) learning techniques, which have gained a lot of popularity recently.
These techniques, despite their superior performance, often lack rigorous constraint satisfaction, which is especially important  in safety-critical applications like medical devices and human-robot interaction, see e.g.~\cite{brunke2022safe}.
We also refer to~\cite{amodei2016concrete}
and~\cite{tambon2022certify} for an overview of the challenges employing learning-based approaches to safety-critical systems; 
and for challenges and recent results in the field of continual learning we refer to the comprehensive surveys~\cite{shaheen2022continual,wang2024comprehensive}.

To address the challenge of ensuring constraint satisfaction while leveraging the benefits of
learning-based control, the field of safe learning has gained prominence and
several safety
frameworks have been proposed~\cite{HewingWaber20, garcia2015comprehensive}, employing various
approaches like control barrier functions~\cite{ames2019control}, Hamilton-Jacobi reachability
analysis~\cite{chen2018hamilton}, Model Predictive Control (MPC)~\cite{aswani2013provably}, and Lyapunov
stability~\cite{perkins2002lyapunov}. Predictive safety filters, as exemplified
in~\cite{Wabersich21,Wabersich23}, verify control input signals against a model to ensure compliance
with prescribed constraints. In~\cite{lanza2023learningbased},
a feedback controller is proposed to compensate for model inaccuracies.
A key feature is that the model can be updated (or even replaced) at runtime while being employed in 
an MPC algorithm.

Due to the activation threshold incorporated in the controller proposed in the current paper,
it may also be used in the context of learning.
When acting as a safety filter for a data-driven learning
algorithm, our controller temporarily interrupts the learning process 
when the activation threshold is surpassed,
resorting to the pure feedback control with ZoH component.
The versatility of our proposed framework is showcased through its
application to prominent data-driven predictive control schemes, specifically data-driven MPC and Reinforcement Learning (RL). 

The data-driven MPC scheme builds upon Willems et al.'s so-called fundamental lemma~\cite{WRMDM05}, allowing a
non-parametric description of the system's input-output behavior based on measurement data, see also~\cite{MarkDorf21,faulwasser2023behavioral} and the references therein. The fundamental lemma states that, for discrete-time linear time-invariant controllable systems, the input-output
trajectories of finite length lie in the column-space of a suitable Hankel matrix constructed
directly from measured input-output data. This result paved the way in the development of data-driven MPC schemes, where the prior model is replaced by measured data,
cf.~\cite{yang2015data,coulson2019data,Berberich20}. 
Therefore, the fundamental lemma is subject to recent substantial research in the field of data-driven control.
In~\cite{SchmitzFaulwasserWorthmann22,faulwasser2023behavioral,tudo:pan22a}, it was extended to
stochastic descriptor systems. Extensions towards continuous-time and non-linear systems were
discussed, e.g., in \cite{lopez2022continuous,rapisarda4370211fundamental,berberich2022linear}
and~\cite{Alsalti21,PersTesi23}, respectively. 

Reinforcement learning has proven to be a successful technique for solving complex 
control problems, e.g.\ single- and multi-agent games~\cite{mnih2015human},
robotics~\cite{ibarz2021train}, and autonomous vehicles~\cite{wang2021data}. The control objective
is usually to either reach a target system state or to maximize the cumulative expected reward,
similar to solving an optimal control problem. Through applying trial-and-error control actions to
the system while collecting data and information during the closed-loop system operation, RL
techniques are able to find a control policy to achieve the desired control task without prior
system knowledge. The main difficulty here is to overcome the exploration-exploitation trade off and to guarantee safety in exploration. 
A comprehensive survey on
applying RL to control systems can be found in~\cite{recht2019tour}. See also the
textbook~\cite{sutton2018reinforcement}~for an overview of Reinforcement Learning, and for its
relationship to optimal control see~\cite{bertsekas2019reinforcement}. 

\noindent
The present article is organized as follows.
In \Cref{Sec:ProblemFormulation}, we provide the problem formulation and 
introduce the system class. 
In \Cref{Sec:MainResult} we introduce the feedback controller component, derive an explicit upper bound on the sampling time~$\tau > 0$, and provide and rigorously proved feasibility result for the ZoH feedback law.
Motivated by a numerical simulation presented in \Cref{Sec:Example}, we combine the proposed feedback ZoH controller with learning-based predictive control algorithms in \Cref{Sec:Extensions}, namely data-driven MPC based on Willems et al.'s 
fundamental lemma in \Cref{Sec:WillemsMPC}, and Reinforcement Learning-based control in \Cref{Sec:RL_Control}.
The integration of the proposed controller into learning-based controllers illustrates its capability to serve as a safety filter for safe online learning.
We prove feasibility of the combined controllers, and demonstrate the superior control performance via numerical simulations.
The more involved proofs, including the proofs of our main results Theorems~\ref{Thm:recursive_feasible_r} and~\ref{Thm:Combined_Controller}, 
are relegated to~\Cref{Sec:Appendix} to make the results more accessible.

\ \\
\textbf{Notation}:
\begin{small}%
$\N,\R$ {is} the set of natural and real numbers, resp.
$\Rp:=[0,\infty)$. 
The standard inner product on $\R^n$ is denoted by $\al \cdot, \cdot \ar$, and $\Norm{x}:=\sqrt{\al x,x\ar}$ {for} $x\in\R^n$.
$\cB_\rho := \setdef{ x \in \R^n }{ \| x \| < \rho }$.
$\con^p(V,\R^n)$ is the linear space of $p$-times continuously  differentiable
functions $f:V\to\R^n$, where $V\subset\R^m$ and $p\in\N\cup \{\infty\}$; $\con(V,\R^n):=\con^0(V,\R^n)$.
For an interval $I\subset\R$,  $L^\infty(I,\R^n)$ {is} the space of measurable essentially bounded
functions $f: I\to\R^n$ with norm $\SNorm{f}:=\esssup_{t\in I}\Norm{f(t)}$. 
$L^\infty_{\text{loc}}(I,\R^n)$ is the space of locally bounded measurable functions.
$W^{k,\infty}(I,\R^n)$ is the Sobolev space of all $k$-times weakly differentiable functions
$f:I\to\R^n$ {with} $f,\dots, f^{(k)}\in L^{\infty}(I,\R^n)$,
{$\Lip(\Rp,\R^m)$ is the space of Lipschitz continuous functions $f: \Rp \to \R^m$.
For a finite sequence $(f_k)_{k=0}^{N-1}$ in $\R^n$ of length $N$ we define the vectorization $f_{[0,N-1]} := \begin{bmatrix}
    f_0^\top & \dots & f_{N-1}^\top
\end{bmatrix}{}^\top \in \R^{nN}$. 
}
\end{small}

\section{Control objective, system class, and preliminary results} \label{Sec:ProblemFormulation}
We consider nonlinear continuous-time control systems
\begin{equation} \label{eq:Sys_r}
\begin{aligned}
y^{(r)}(t) &= f \big( d(t), \oT(y,\ldots,y^{(r-1)} )(t) \big) \\
& \qquad \quad + g \big( d(t), \oT(y,\ldots,y^{(r-1)} )(t) \big) u(t), \\
y|_{[-\sigma,0]} &= y^0  \in \cC^{r-1}([-\sigma,0],\R^m),
\end{aligned}
\end{equation}
where $d \in L^\infty (\Rp, \R^p)$ represents an unknown bounded disturbance, $f \in \cC(\R^p \times \R^q, \R^m)$ is a drift term, the function $g \in \cC(\R^p \times \R^q, \R^{m \times m})$ is the input gain function, and the operator $\oT$ is causal,
locally Lipschitz and satisfies a bounded-input bounded-output property; the operator is characterized in detail in \Cref{Def:Operator-class}, and the class of systems under consideration is introduced in \Cref{Def:system-class}.
We emphasize that many physical phenomena such as \emph{backlash} and \emph{relay hysteresis}, and \emph{nonlinear time delays} can be modeled by means of the operator~$\oT$ ($\sigma$ corresponds to the initial delay), {cf.~\cite[Sec.~1.2]{BergIlch21}.}
Moreover, systems with infinite-dimensional internal dynamics can be represented by~\eqref{eq:Sys_r}, see e.g. \cite{BergPuch20}.
For a control function $u\in L^\infty_{\loc}(\Rp, \R^m)$, system~\eqref{eq:Sys_r} has a solution in 
the sense of \textit{Carath\'{e}odory}, meaning a function
${x} : [-\sigma,\omega) \to \R^{rm}$, $\omega > 0$, 
with ${x}|_{[-\sigma,0]} = (y^0,\dot y^0,\ldots,(y^0)^{(r-1)})$ 
such that ${x}\vert_{[0,\omega)}$ is absolutely continuous and 
satisfies ${\dot x}_i(t) = {x}_{i+1}(t)$ for $i=1,\ldots,r-2$, and ${\dot x}_r(t) = f(d(t),\oT({x}(t))) + g(d(t),\oT({x}(t))) u(t)$ (which corresponds to~\eqref{eq:Sys_r} with $y=x_1$) for almost all~$t\in[0,\omega)$.
A solution ${x}$ is said to be \textit{maximal}, if it does not have a right extension which is also a solution.

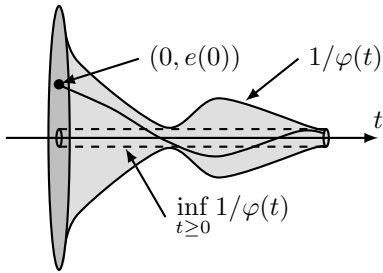
\begin{figure}[hbt]
  \begin{center}
\begin{tikzpicture}[scale=0.35]
\tikzset{>=latex}
  \filldraw[color=gray!25] plot[smooth] coordinates {(0.15,4.7)(0.7,2.9)(4,0.4)(6,1.5)(9.5,0.4)(10,0.333)(10.01,0.331)(10.041,0.3) (10.041,-0.3)(10.01,-0.331)(10,-0.333)(9.5,-0.4)(6,-1.5)(4,-0.4)(0.7,-2.9)(0.15,-4.7)};
  \draw[thick] plot[smooth] coordinates {(0.15,4.7)(0.7,2.9)(4,0.4)(6,1.5)(9.5,0.4)(10,0.333)(10.01,0.331)(10.041,0.3)};
  \draw[thick] plot[smooth] coordinates {(10.041,-0.3)(10.01,-0.331)(10,-0.333)(9.5,-0.4)(6,-1.5)(4,-0.4)(0.7,-2.9)(0.15,-4.7)};
  \draw[thick,fill=lightgray] (0,0) ellipse (0.4 and 5);
  \draw[thick] (0,0) ellipse (0.1 and 0.333);
  \draw[thick,fill=gray!25] (10.041,0) ellipse (0.1 and 0.333);
  \draw[thick] plot[smooth] coordinates {(0,2)(2,1.1)(4,-0.1)(6,-0.7)(9,0.25)(10,0.15)};
  \draw[thick,->] (-2,0)--(12,0) node[right,above]{\normalsize$t$};
  \draw[thick,dashed](0,0.333)--(10,0.333);
  \draw[thick,dashed](0,-0.333)--(10,-0.333);
  \node [black] at (0,2) {\textbullet};
  \draw[->,thick](4,-3)node[right]{\normalsize$\inf\limits_{t \ge 0} 1/\vp(t)$}--(2.5,-0.4);
  \draw[->,thick](3,3)node[right]{\normalsize$(0,e(0))$}--(0.07,2.07);
  \draw[->,thick](9,3)node[right]{\normalsize $1/\vp(t)$}--(7,1.4);
\end{tikzpicture}
\end{center}
 \vspace*{-2mm}
 \caption{Error evolution in a funnel $\mathcal F_{\vp}$ with boundary $1/\vp(t)$; the figure is based on~\cite[Fig.~1]{BergLe18a}, edited for present purpose.}
 \label{Fig:funnel}
 \vspace{-4mm}
 \end{figure}

\noindent
The control objective is to design a
zero-order hold control strategy, i.e., for sampling time~$\tau > 0$, 
\[
    u(t) \equiv u \qquad \forall\, t \in [t_i, t_i + \tau),
    \ i \in \N,
\]
where the data are collected at uniform sample times ${t_i = i \cdot \tau \in \Rp}$,
which achieves for a system~\eqref{eq:Sys_r} output tracking of a given reference
$y_{\rf}\in W^{r,\infty}(\Rp,\R^{m})$ 
within pre-specified error bounds. 
To be more
precise, the tracking error $t\mapsto e(t):=y(t)-y_{\rf}(t)$ shall evolve within the prescribed
performance funnel
\begin{align*}
    \cF_\phi= \setdef{(t,e)\in \Rp\times\R^{m}}{\phi(t)\Norm{e} < 1}.
\end{align*}
This funnel is determined by 
the function~$\vp$ belonging to
\[
     \cG:=
     \setdef
         {\vp\in W^{1,\infty}(\Rp,\R)}
         {
            \inf_{s \ge 0} \vp(s) > 0
         },
 \]
see \Cref{Fig:funnel} for an illustration.

The specific application usually dictates the constraints on the tracking error and thus
indicates suitable choices for~$\vp$.
To achieve the control objective, we introduce auxiliary error variables.
For $\phi\in\cG$, ${y_{\rf}\in W^{r,\infty}(\Rp,\R^m)}$,
a bijection $\alpha \in\con^1([0,1),[1,\infty))$,
$t\geq0$, and ${{\xi}:=({\xi}_1,\ldots,{\xi}_r)\in\R^{rm}}$ we formally introduce the following auxiliary variables 
\begin{align} \label{eq:ek}
    e_1(t,{\xi})&:=\phi(t)({\xi}_1-y_{\rf}(t)),\\
    e_{k+1}(t,{\xi})&:=\phi(t)({\xi}_{k+1}-y_{\rf}^{(k)}(t))\!+\!\alpha(\Norm{e_{k}(t,{\xi})}^2)e_{k}(t,{\xi}), \nonumber
    \end{align}
for $k=1,\ldots,r-1$,
where $e_1(t)$ is the tracking error~$e(t)$ normalized with respect to the error boundary~$\vp(t)$.
A suitable choice for the bijection is $\alpha(s):=1\slash(1-s)$.
While the error $e_{k+1}$ is formally not defined for $\Norm{e_{k}(t,{\xi})}=1$, 
we will in the following only evaluate the error variables for $\xi\in\cD_{t}^r$ as defined in~\eqref{eq:DefSetD} in \Cref{Sec:Aux_results} excluding this edge case. In favor of a simpler notation, we therefore refrain from defining $e_{k+1}$ at these points. 
Using the short notation $e_r(t):=e_r(t,(y,\dot{y},\ldots,y^{(r-1)})(t))$,
we propose the following controller structure for~${i \in \N}$
    \begin{equation} \label{eq:controller_recursive}
   \forall  t \in  [t_i, t_i + \tau)  : \, u(t) =
    \begin{cases}
     0, \!\! &  \|  e_r(t_i)\| < \lambda, \\
         - \beta \tfrac{e_r(t_i)}{\|e_r(t_i)\|^2}  , \!\! &\|  e_r(t_i)\| \ge \lambda,
    \end{cases}
    \end{equation}
where $\lambda \in (0,1)$ is an {activation threshold}, and~$\beta > 0$ is the input gain.
{In \Cref{Sec:Aux_results} we show $ e_r \in \cB_1$. Thus,} the control function~$u$ is uniformly bounded since we have
\begin{align*}
   \forall t \geq 0:\quad \|u(t) \| \leq \frac{\beta}{\lambda}.
\end{align*}
Our designed controller can be considered to be similar to funnel control, see \cite{BergIlch21,BergLe18a,IlchRyan02b}, in terms of its ability to achieve output reference tracking within predefined error boundaries, as well as concerning the used intermediate error variables~\eqref{eq:ek}.
On the other hand, contrary to the standard funnel controller, the feedback law~\eqref{eq:controller_recursive} is a normalized linear sample-and-hold output feedback with uniform sampling rate.
Since it involves an activation threshold, it has also similarity with the zero-or-hold controller in~\cite{Schenato09}.
A further essential difference to continuous funnel control is that in the present approach the control objective is achieved by using estimates about the system dynamics, while in funnel control no such information is used to the price that the maximal control effort cannot be estimated a-priori.

\subsection{System class} \label{Sec:SystemClass}
In this section we formally introduce the system class under consideration.
Prior to that, we state assumptions on the system parameters and characterize the operator~$\oT$.
\begin{assumption} \label{Ass:d_bounded}
    A bound $D>0$ for the \emph{unknown} disturbance
    $d \in L^\infty(\Rp, \R^p)$ with $\|d\|_\infty \le D$ is known.
\end{assumption}
\begin{assumption} \label{Ass:g_sign_definite}
The matrix valued function {$g\in \cC(\R^p \times \R^q, \R^{m \times m})$} is strictly positive definite, that is
\[
    \forall {x} \in \R^{p + q} \ \forall z \in \R^m \setminus \{0\} \, : \  \al z, g({x}) z \ar > 0.
\]
\end{assumption}
Note that we could also allow the case of strictly negative $g$ by changing the sign in \eqref{eq:controller_recursive}.
Further, note that, while some authors only use the term \textit{strictly positive definite} for symmetric matrices,
we do not assume $g(x)$ to be symmetric.
Next, we provide the defining properties of the class of operators to which $\oT$ in~\eqref{eq:Sys_r} belongs.
\begin{definition} \label{Def:Operator-class} 
For $n,q\in\N$ and $\sigma\geq 0$, the set $\cT^{n,q}_{\sigma}$ denotes the class of operators 
$\textbf{T}:\con([-\sigma,\infty),\R^{n})\to L^\infty_{\loc} (\Rp, \R^{q})$
for which the following properties hold:
\begin{enumerate}[i)]
    \item\textit{Causality}:  $\fa y_1,y_2\in\con([-\sigma,\infty),\R^{n})$  $\fa t\geq 0$:
    \[
        y_1\vert_{[-\sigma,t]} = y_2\vert_{[-\sigma,t]}
        \quad \Impl\quad
        \textbf{T}(y_1)\vert_{[0,t]}=\textbf{T}(y_2)\vert_{[0,t]}.
    \]
    \item\textit{Local Lipschitz}: 
    $\fa t \ge 0 $ $\fa y \in \con([-\sigma,t] ; \R^n)$ 
    $\ex \Delta, \delta, c > 0$ 
    $\fa y_1, y_2 \in \con([-\sigma,\infty) ; \R^n)$ with
    $y_1|_{[-\sigma,t]} = y = y_2|_{[-\sigma,t]}$
    and $\Norm{y_1(s) - y(t)} < \delta$,  $\Norm{y_2(s) - y(t)} < \delta $ for all $s \in [t,t+\Delta]$:
    \[
        \esssup_{\mathclap{s \in [t,t+\Delta]}}  \Norm{\textbf{T}(y_1)(s) \!-\! \textbf{T}(y_2)(s) }  
        \!\le\! c \ \sup_{\mathclap{s \in [t,t+\Delta]}} \Norm{y_1(s) - y_2(s)}\!.
    \]
    \item\textit{Bounded-input bounded-output (BIBO)}:
    $\fa c_0 > 0$ $\ex c_1>0$  $\fa y \in \con([-\sigma,\infty), \R^n)$:
    \[
    \sup_{t \in [-\sigma,\infty)} \Norm{y(t)} \le c_0 \ 
    \Impl \ \sup_{t \in [0,\infty)} \Norm{\textbf{T}(y)(t)}  \le c_1.
    \]
\end{enumerate}
\end{definition}
While the first property (causality) introduced in
\Cref{Def:Operator-class} is quite intuitive, the second (locally Lipschitz) is of a more technical
nature, required to guarantee existence and uniqueness of solutions. The third property (BIBO)
can be motivated from a practical point of view as an infinite-dimensional
extension of minimum-phase.
Various examples for the operator~$\oT$ can be found in~\cite{IlchRyan02b,BergIlch21}.

With \Cref{Ass:d_bounded,Ass:g_sign_definite} and \Cref{Def:Operator-class}
we formally introduce the system class under consideration.
\begin{definition} \label{Def:system-class}
     For $m,r \in \N$ a system~\eqref{eq:Sys_r} belongs to the system class $\cN^{m,r}$, written ${(d,f,g,\oT)
     \in\cN^{m,r}}$, if, for some $p,q\in\N$ and $\sigma \geq0$, the following holds:
     ${d \in L^\infty(\Rp,\R^p)}$ satisfies \Cref{Ass:d_bounded},
     ${f\in\con(\R^p\times \R^q ,\R^m})$,
     $g$ satisfies \Cref{Ass:g_sign_definite},
     and ${\oT\in\cT^{rm,q}_{\sigma}}$.
\end{definition}
Note that all linear minimum-phase systems with relative degree $r \in \N$ are contained in the system class~$\cN^{m,r}$, cf.~\cite{BergIlch21}.
Moreover, under assumptions provided in~\cite[Cor.~5.6]{ByrnIsid91a},
   a nonlinear system of the form
   \begin{equation}\label{eq:SysWithState}
   \begin{aligned}
       \dot{x}(t)  & = \tilde f(x(t)) + \tilde g(x(t)) u(t),\quad x(0)=x^0\in\R^n,\\
       y(t)        & = h(x(t)),
   \end{aligned} 
   \end{equation}
   with  nonlinear functions~$\tilde f:\R^n\to \R^n$,
   $\tilde g:\R^n\to \R^{n\times m}$ and $h : \R^n \to \R^m$, 
   can be put in the form~\eqref{eq:Sys_r} with $\sigma=0$ and appropriate functions~$f$ and~$g$ and an operator $\textbf{T}$
   via  a coordinate transformation induced by a diffeomorphism~$\Phi:\R^n\to\R^n$.
   The  operator~$\oT$ then is the solution operator of the internal dynamics of the transformed system.
   Using the diffeomorphism~$\Phi$, the presented results this paper can also be expressed for the system~\eqref{eq:SysWithState}.
   In this case, exact knowledge about the coordinate transformation is not required for 
   the design and application of the presented controller -- it merely serves as a tool for the proofs. 

\subsection{Auxiliary results} \label{Sec:Aux_results}
In order to formulate the main result about feasibility of the proposed ZoH controller,  we introduce some notation and establish two auxiliary results in this section.
{
We use the shorthand notation 
\begin{equation*}
    \chi(y)(t):=(y(t),\dot{y}(t),\ldots,y^{(r-1)}(t))\in\R^{rm}    
\end{equation*}
for $y\in W^{r,\infty}(\Rp,\R^m)$ and $t\in\Rp$. 
}
To guarantee that the tracking error $e = y - y_{\rm \rf}$ evolves within the boundary of~$\cF_\phi$, we want to address the problem 
of ensuring that $\chi(y)(t)$ is at every time $t\geq 0$ an element of the set
\begin{equation}\label{eq:DefSetD}
        \cD_{t}^{r}:=\setdef
        {\!\! {\xi}\in\R^{rm}\!\!}
        { \!\!\!\!\begin{array}{l}
              \Norm{e_{k}(t,{\xi})}\!<1,\ k=1,\ldots, r-1, \\
              \|e_r(t,{\xi})\| \le 1
        \end{array} \!\!\!\!} \!.
\end{equation}
We define the set of all functions $\zeta\in\con^r([-\sigma,\infty),\R^m)$
which coincide with~$y^0$ on the interval $[-\sigma,0]$ and for which 
$\chi(y)(t)\in\cD^r_t$ on the interval $[t_0,\delta)$ for $\delta\in(0,\infty]$:
\[
    \cY^r_\delta\!\!:=\!\setdef
        {\!\!\!\zeta\in \con^{{r-1}}\!([-\sigma,\!\infty),\R^m)\!\!\!}
        {\!\!\!\!\begin{array}{l}
             \zeta|_{[-\sigma,0]}=y^0,  \\
             \!\!\fa t\in [0,\delta):\chi(\zeta)(t)\in\cD_t^r\!\!
        \end{array}\!\!\!\!}\!.
\]
We aim to infer the existence of bounds for the error variables $e_k$ defined in~\eqref{eq:ek} for all functions in~$\cY_{\delta}^r$ independent of the {functions $f$, $g$, the disturbance $d$, the operator $\oT$, and the applied control $u$ in} system dynamics~\eqref{eq:Sys_r}.
To this end,
we introduce the following constants
$\ve_k,\mu_k$.
Let $\ve_0=0$ and $\bar\gamma_0 :=0$.
Successively for $k={1},\ldots,r-1$ define 
\begin{align}
\hat \ve_k \!&\in\! (0,1)  \!\text{ \rm s.t. } 
 \alpha(\hat \ve_k^2) \hat \ve_k \!=\!  \SNorm{\frac{\dot \vp}{\vp} }\!\!\!\!\!\! ( 1\! +\! \alpha(\ve_{k-1}^2) \ve_{k-1})\! +\! 1 \!+\! \bar \gamma_{k-1}, \nonumber \\
    \ve_k \!&:= \max \{ \| e_k(0)\|,  \hat \ve_k\} < 1, \label{eq:ve_mu_gam} \\
    \mu_k \!& := \SNorm{\frac{\dot \vp}{\vp} }\!\!\! ( 1 \! +\!  \alpha(\ve_{k-1}^2) \ve_{k-1} ) \! + \! 1\!+\! \alpha(\ve_k^2) \ve_k  \! + \! \bar \gamma_{k-1} , \nonumber \\
    \bar \gamma_k & := 2 \alpha'(\ve_k^2) \ve_k^2 \mu_k + \alpha(\ve_k^2) \mu_k. \nonumber
\end{align}
With these constants we may derive the following result.
\begin{lemma} \label{Lemma:e_k}
Let {$y_{\rf}\in W^{r,\infty}(\Rp,\R^m)$}, $\vp \in \cG$, and $y^0\in\con^{r-1}([-\sigma,0],\R^m)$ with $\chi(y^0)\in\cD_0^r$ be given. 
Then there exist constants $\ve_k,\mu_k>0$ defined in~\eqref{eq:ve_mu_gam} such that for all $\delta\in(0,\infty]$ and  all $\zeta\in \cY^r_\delta$ 
the functions $e_k$ defined in~\eqref{eq:ek} satisfy
\begin{enumerate}[i)]
    \item $\| e_k(t,\chi(\zeta)(t))\|  \ {\leq}\ \ve_k < 1$,
    \item $\| \dd{t} e_k(t,\chi(\zeta)(t)) \| \ {\leq}\ \mu_k$,
\end{enumerate}
for all $t\in [0,\delta)$ and for all $k=1,\ldots,r-1$.
\end{lemma}
The proof is relegated to the \Cref{Sec:Appendix}.
Next, we derive bounds on the right-hand side of system~\eqref{eq:Sys_r}.
\begin{lemma}\label{Lemma:DynamicsBounded}
Consider~\eqref{eq:Sys_r} with $(d,f,g,\oT)\in\cN^{m,r}$.
Let {$y_{\rf} \in W^{r,\infty}(\Rp,\R^m)$}, $\phi\in\cG$, $y^0\in\con^{r-1}([-\sigma, 0],\R^m)$ with $\chi(y^0)(0)\in\cD_{0}^r$, and $D>0$ from \Cref{Ass:d_bounded}.
Then, there exist constants $\fM$,  $\gM$,  $\gm>0$ such that for every $\delta\in(0,\infty]$,
$\zeta\in\cY^r_{\delta}$, $d \in L^\infty(\Rp,\R^p)$ with $\SNorm{d}\leq D$, $z\in\R^{m}\backslash\cbl0\cbr$, and $t\in[0,\delta)$ 
\begin{equation} \label{eq:fmax_gmax_gmin}
\begin{aligned}
    \fM &\geq \SNorm{f((d,\oT(\chi(\zeta)))|_{[0,\delta)})} ,\\
    \gM &\geq \SNorm{g((d,\oT(\chi(\zeta)))|_{[0,\delta)})} ,\\
    \gm &\leq \frac{\al z, g((d,\oT(\chi(\zeta)))|_{[0,\delta)}(t))z\ar }{\Norm{z}^2}.
\end{aligned}
\end{equation}
\end{lemma}
The proof is relegated to the \Cref{Sec:Appendix}.
\Cref{Lemma:DynamicsBounded} ensures existence of bounds on the dynamics of the system to be controlled.
To compute these bounds, \textit{some} system knowledge is necessary.
For instance, if the structure of the governing equations is known and the parameters are known to be in
a certain range, the worst case estimates $\fM, \gM, \gm$ can be computed using the desired reference trajectory and the prescribed error tolerance, i.e., seeking the maximum of continuous functions within a compact set.

\section{Sampled-data feedback controller} \label{Sec:MainResult}
With the introductory results presented in the previous section, we are now in a position to formulate a feasibility result about the ZoH feedback controller.
To phrase it, \Cref{Thm:recursive_feasible_r} yields that the ZoH controller~\eqref{eq:controller_recursive} achieves the control objective discussed in \Cref{Sec:ProblemFormulation} for a system~\eqref{eq:Sys_r} with $(d,f,g,\oT) \in \cN^{m,r}$, if the sampling time~$\tau$ satisfies the following condition~\eqref{eq:tau_r}.
\begin{theorem} \label{Thm:recursive_feasible_r} 
Given a reference {$y_{\rf} \in W^{r,\infty}(\Rp,\R^m)$} and a funnel function $\vp \in \cG$ consider a system~\eqref{eq:Sys_r} with $(d,f,g,\oT) \in \cN^{m,r}$.
With the constants given in~\eqref{eq:ve_mu_gam}, set
\[
        \kappa_0 \!\!:=\!\! \SNorm{\frac{\dot \vp}{\vp}}\!\!\!\!\!\! ( 1 + \alpha(\ve_{r-1}^2) \ve_{r-1} )
                    + \SNorm{\vp}\!\! ( \fM + \|y_{\rf}^{(r)}\|_\infty )+ \bar \gamma_{r-1},
 \]
define the input gain
    \begin{equation*}
        \beta > \frac{2 \kappa_0}{ \gm \inf_{s \ge 0} \vp(s) },
    \end{equation*}
and the constant
$
   \kappa_1 := \kappa_0 +  \SNorm{\vp} \gM \beta.
$
Assume that the initial condition satisfies
    $\chi(y^0)(0) \in \cD_{0}^r$, i.e., the error variables from~\eqref{eq:ek} (here we omit the dependence on $\chi(y) = (y,\ldots,y^{(r-1)})$) satisfy
    $\| e_k(0)\| < 1$ for all $k=1,\ldots,r-1$, and $e_r(0) \le 1$;
and, for an activation threshold ${\lambda \in (0,1)}$, let the sampling time satisfy
\footnote{
In the published version \url{https://doi.org/10.1016/j.sysconle.2024.105892}, there is a typo: $\lambda^2$ is missing in~\eqref{eq:tau_r}.
}
    \begin{equation} \label{eq:tau_r}
0 <      \tau \le \min \left\{ \frac{ {\lambda^2} \kappa_0 }{\kappa_1^2}, \frac{1-\lambda}{\kappa_0} \right\}.
    \end{equation}
Then  the ZoH controller~\eqref{eq:controller_recursive}
applied to a system~\eqref{eq:Sys_r} yields that $\| e_k(t)\| < 1$ for all $k=1,\ldots,r-1$ and $\|e_r(t)\| \le 1$ for all $t \ge 0$.
This is initial and recursive feasibility of the ZoH control law~\eqref{eq:controller_recursive}.
In particular, the tracking error satisfies $\|e(t)\| < 1/\vp(t)$ for all $t \in \Rp$.
\end{theorem}
The proof of \Cref{Thm:recursive_feasible_r} is relegated to the \Cref{Sec:Appendix}.

    The parameter $\lambda \in (0,1)$ in~\eqref{eq:controller_recursive} is an activation threshold (cf. event-triggered control~\cite{Heemels2021}), chosen by the designer, which divides the tracking error in a safe and a safety critical region.
    A large value of $\lambda$ implies that the controller will be inactive for a wide range of values of the last error variable, which, in case of relative degree one, means inactivity for a wide range of the tracking error, while still guaranteeing transient accuracy.

The sampling time $\tau$ in~\eqref{eq:tau_r} strongly depends on the {evolution of the} funnel function and on the reference~$y_{\rf}$.
     This gives the possibility of dynamically adapting the sampling time, e.g., in the case of setpoint transition, where the reference is constant $y_{\rf}^0$ in the first period and constant $y_{\rf}^1 \neq y_{\rf}^0$ in the last period. 
    At the setpoints the sampling time can be larger than during the transition.

    An explicit bound on the control input can be computed in advance, since $\| u \|_\infty \le \beta/\lambda$. 
    This bound depends on the system parameters derived in \Cref{Lemma:DynamicsBounded}.
    However, precise knowledge about the functions $f$, $g$ and the operator~$\oT$ is not necessary.
    Mere (conservative) estimates on the bounds $\fM$, $\gM$, and $\gm$ as in~\eqref{eq:fmax_gmax_gmin} are sufficient.

\begin{remark}
    The results in \Cref{Thm:recursive_feasible_r} are also valid for $\|e_r(0)\| = 1$.
    This is in contrast to continuous time funnel control, where all $r$ error variables~\eqref{eq:ek} initially have to be bounded away from~$1$ to guarantee boundedness of the input.
    To illustrate this, consider $\dot y(t) = u(t)$, and~$y_{\rm ref} = 0$.
    Let $\vp \in \cG$ and choose the bijection $\alpha(s) = 1/(1-s)$. 
    According to \cite{BergIlch21} the control is given by $u(t) = - \tfrac{y(t)}{1-\vp(t)^2y(t)^2}$.
    Now, for a sequence of initial values $y_j(0)$, $j \in \N$, such that $\vp(0) |y_j(0)| \rightarrow 1$ for $j \to \infty$, the sequence of corresponding initial controls $u_j(0)$ is unbounded.
    On the other hand, for the same  sequence of initial values the controller~\eqref{eq:controller_recursive} yields a bounded signal $\|u_{\rm ZoH}\|_\infty \le \beta / \lambda$.
    Moreover, such a sequence of initial values requires ever smaller sampling time, if a continuous funnel controller is implemented, cf.~\Cref{Sec:Example}.
\end{remark}

    \begin{remark} \label{Rem:u_eq_0_not_nec}
    Note that $u=0$ is not necessary for $\|e_r(t_i)\| < \lambda$; however,
    according to the current proof, $u \neq 0$ will decrease~$\tau$.
    For instance, applying the control value~$u(t_{i-1})$ of the last sampling period is feasible, or the control value may be chosen according to the data informativity framework~\cite{van2020data}.
    Such a data-driven control is safeguarded by the proposed controller~\eqref{eq:controller_recursive}, similar to the combined controller~\cite{lanza2023learningbased}.
            We will exploit this observation in \Cref{Sec:Extensions}, where we propose a two-component data-driven/learning-based controller with ${u \neq 0}$ for ${\|e_r(t_k)\| < \lambda}$.
    \end{remark}

\section{Numerical example: pure ZoH feedback} \label{Sec:Example}
\begin{figure}[h!]
\begin{center}
\includegraphics[trim=2cm 4cm 5cm 15cm,clip=true,width=4.3cm]{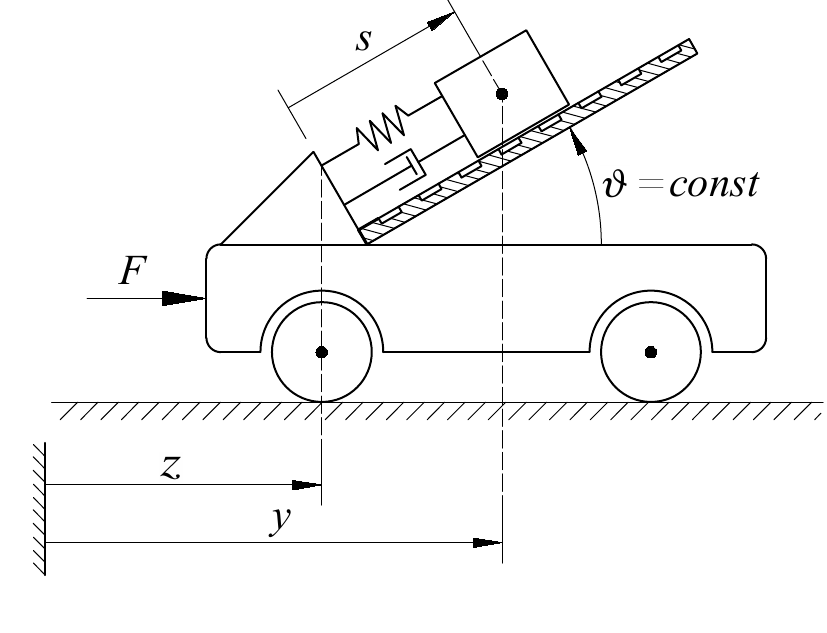}
\end{center}
    \vspace{-5mm}
    \caption{Mass-on-car system. The figure is based on~\cite{SeifBlaj13,BergIlch21}.}
    \label{Fig:Mass-on-a-car}
    \vspace{-2mm}
\end{figure}

To illustrate the controller~\eqref{eq:controller_recursive} we consider the mass-on-car system~\cite{SeifBlaj13}.
On a car with mass~$m_1$, to which a force~$F=u$ can be applied, a ramp is mounted on which a second mass~$m_2$ moves passively, see \Cref{Fig:Mass-on-a-car}.

\begin{figure}
    \centering
    \includegraphics[scale=0.32]{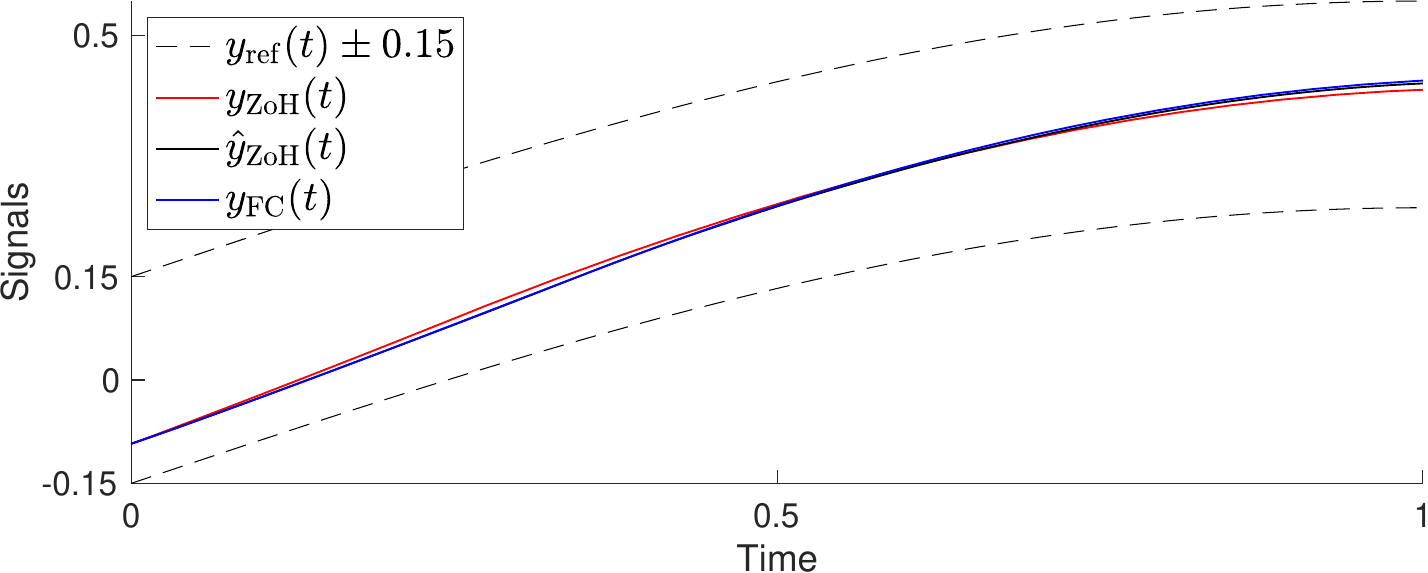}
    \vspace{-2mm}
    \caption{Outputs, reference, and error tolerance.}
    \label{fig:ZOH_Errors}
    \vspace{-4mm}
\end{figure}

\begin{figure}
    \centering
    \includegraphics[scale=0.32]{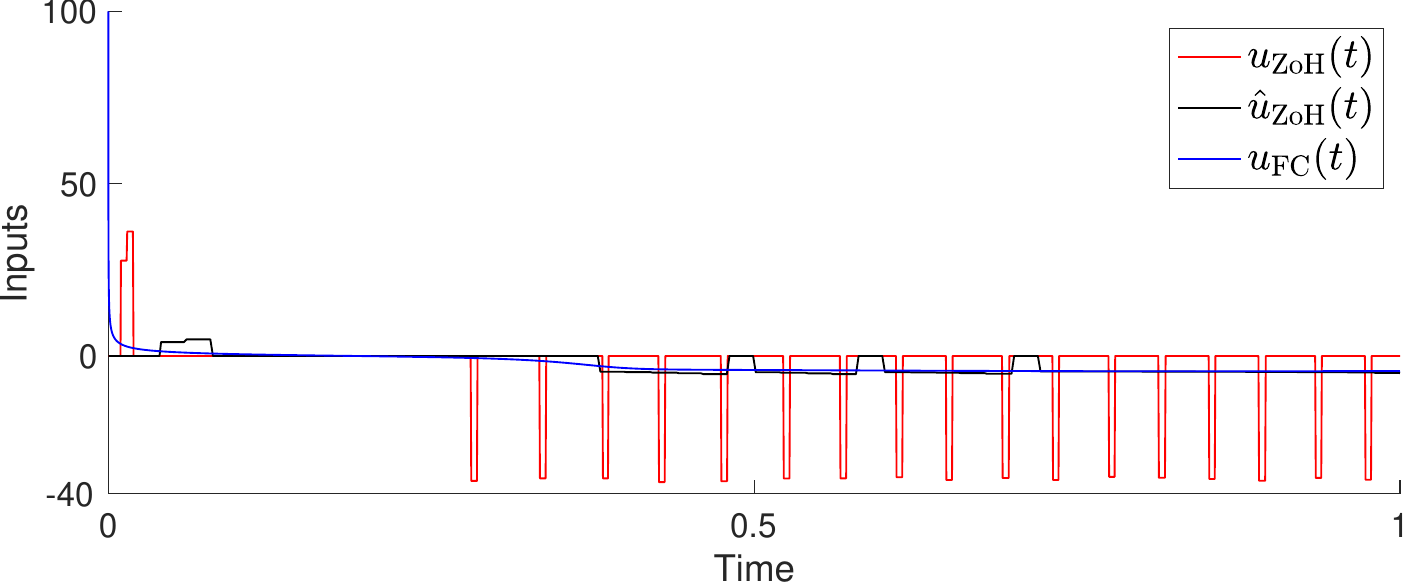}
    \vspace{-2mm}
    \caption{Controls.}
    \label{fig:ZOH_Controls}
    \vspace{-6mm}
\end{figure}

The second mass is coupled to the car by a spring-damper combination, and the ramp is inclined by a fixed angle~$\vartheta \in (0,\pi/2)$.
The equations of motion are given by
\begin{subequations} \label{eq:MOC}
\begin{equation}
    \begin{aligned}
    \resizebox{0.88\hsize}{!}{$
        \begin{bmatrix}
        m_1 + m_2 & m_2 \cos(\vartheta) \\
        m_2 \cos(\vartheta) & m_2
        \end{bmatrix}
        \! \!
        \begin{pmatrix}
        \ddot z(t) \\ \ddot s(t)
        \end{pmatrix}
        \!+\!
        \begin{pmatrix}
            0 \\
            k s(t) + d \dot s(t) 
        \end{pmatrix}
        = 
        \begin{pmatrix}
            u(t) \\ 0
        \end{pmatrix},
        $}
    \end{aligned}
\end{equation}
where~$z$ is the car's horizontal position, and~$s$ is the relative position of the second mass.
As output the second mass' horizontal position is measured
\begin{equation}
    y(t) = z(t) + \cos(\vartheta) s(t).
\end{equation}
\end{subequations}
For simulation we choose the parameters~$\vartheta = \pi/4$, $m_1 = 1$, $m_2 = 2$, spring constant $k=1$, and damping~$d=1$.
A short calculation yields that for these parameters system~\eqref{eq:MOC} has relative degree $r=2$, and
as outlined in~\cite[Sec.~3]{BergIlch21} it can be represented in the form~\eqref{eq:Sys_r} with BIBO internal dynamics.
We simulate output reference tracking of the signal $y_{\rf} = 0.4 \sin(\tfrac{\pi}{2}t)$ for $t \in [0,1]$, transporting the mass~$m_2$ on the car from position~$0$ to $0.4$ within chosen error boundaries~$\pm 0.15$.
We choose the activation threshold~${\lambda = 0.75}$.
With these parameters a brief calculation (using the variation of constants formula for the internal dynamics) yields $\fM \le 1.4$, $\gm = \gM = 0.25$, and hence, the sampling time $\tau \le 4.8 \cdot 10^{-3}$, and the gain $\beta \ge 27.55$, which guarantee success of the tracking task. Choosing the smallest~$\beta$ this already gives $\|u_{\rm ZoH}\|_\infty \le \beta/\lambda  \le 36.73$.
We start with a small initial tracking error $y(0) = -0.0925$, and $\dot y(0) = \dot y_{\rm ref}(0)$.
We compare the controller~\eqref{eq:controller_recursive} with the continuous funnel controller~\cite{BergIlch21}; corresponding signals have the subscript $\rm FC$, e.g., $u_{\rm FC}$.
Moreover, simulating the ZoH controller was even successful for $\tau = 2.0\cdot 10^{-2}$ and $\beta = 4$; corresponding signals have a circumflex, e.g., $\hat y_{\rm ZoH}$.
\Cref{fig:ZOH_Errors} shows the system's output and the reference plus/minus error tolerance.
Note that although the control input is discontinuous, the output signal is continuous due to integration.
All controllers achieve the tracking task.
In \Cref{fig:ZOH_Controls} the controls are depicted.
The ZoH input consists of separated pulses - for two reasons. 
First, the control law~\eqref{eq:controller_recursive} uses (undirected) worst-case estimations $\gm, \gM$ and $\fM$ to compute the input signal. Hence, the control signal is at many time instances unnecessary large; however, it is ensured that the control signal is sufficiently large for all times.
Second,~\eqref{eq:controller_recursive} involves an activation threshold~$\lambda$, i.e., the controller is inactive, if the tracking error is small.
If at sampling the tracking error is above this threshold, the applied input is sufficiently large (due to the worst case estimations) to push the error back below the threshold. Thus, at the next sampling instance the input is determined to be zero.
The worst-case estimations and the ZoH setting make it inevitable that the control signal looks peaky. 
The control signal $\hat u_{\rm ZoH}$ (black) is also peaky, but not so large in magnitude (smaller~$\beta$) and with a larger width (larger~$\tau$). Overall, $\hat u_{\rm ZoH}$ is comparable with $ u_{\rm FC}$.
The success of the simulation with these parameters gives rise to the hope of finding better estimates for sufficient control parameters~$\beta,\tau$ in future work.
Improving the control performance is also topic of \Cref{Sec:Extensions}.
Note that the control signal~$u_{\rm FC}$ also has a large peak at the beginning, where $\|u_{\rm FC}\|_\infty \approx 100$.
For simulation, we used \textsc{Matlab}, for integration of the dynamics the routine \texttt{ode15s} with $\rm AbsTol=RelTol=10^{-6}$, with adaptive step size.
Integrating the funnel controller~\cite{BergIlch21} \texttt{ode15s} yields that the maximal step size is $\approx 3.99 \cdot 10^{-2}$ and the minimal step size is $\approx 1.21 \cdot 10^{-6}$.
This means, the largest step is about eight times larger than $\tau$, and the smallest time step is about 4000 times smaller than~$\tau$.

\section{Safeguarded data-based control} \label{Sec:Extensions}
As can be seen from the numerical simulation in \Cref{Sec:Example}, the control signal $u_{\rm ZoH}$ exhibits undesirably large peaks.
This is due to the worst case estimations in the controller design.
In this section, a basic idea for improving the control signal is explained using two example techniques.

\begin{figure}[H]
  \begin{center}
\begin{tikzpicture} [auto, node distance=2cm,>=stealth, every text node part/.style={align=center}]
\def\hoch{0.8cm};
\def\breit{0cm};
\def\distu{1.3cm};
\def\dista{1.3cm};
\node [block, minimum width = \breit, minimum height = \hoch,] (System) {System~\eqref{eq:Sys_r}};
\node [block, minimum width = 3.8cm, minimum height = \hoch, below of=System, node distance = 1.3cm, xshift=-2mm ] (ZoH) { $u_{\rm ZoH}$ from~\eqref{eq:controller_recursive} \\ (\emph{safety critical region})};
\node [block, minimum width = 3.8cm, minimum height = \hoch, below of=ZoH, node distance = 1.2cm ] (ZoH+MPC) {Learning-based control \\ (\emph{safe region})};
\node[input, left of=System, node distance = 2.3cm, xshift=-2mm] (u) {};
\coordinate[right of=u, node distance = 0.0cm] (uin) {};
\coordinate[below of=uin, node distance = 1.9cm] (uin_control) ;
\node[circle,,draw=black, fill=white,inner sep=0pt,minimum size=3pt,  right of=System, node distance = 3.8cm] (ref_in) {$-$};
\node[right of=ref_in, node distance = 1.8cm] (ref) {};
\coordinate[below of=ref_in, node distance = 1.9cm] (er);
\draw[->] (u) --  node[above]{$u$} (System);
\draw[->] (System) --  node[above]{$y$} (ref_in);
\draw[->] (ref) --  node[above]{$y_{\rm ref}$} (ref_in);
\draw[->] (er) --  node[yshift=3pt,above]{$\|e_r(t_k)\| \ge \lambda$} (ZoH);
\draw[->] (er) --  node[yshift=-4pt,below]{$\|e_r(t_k)\| < \lambda$} (ZoH+MPC);
\draw[-] (ref_in) -- node[right]{$e = y - y_{\rm ref}$,\\ $e_k$ as in~\eqref{eq:ek}} (er);
\draw[-] (ZoH) -- (uin_control) --  (uin) ;
\draw[-] (ZoH+MPC) -- (uin_control);
\end{tikzpicture}
\end{center}
 \vspace*{-2mm}
 \caption{Schematic structure of the combined controller.}
 \label{Fig:Controller}
 \vspace*{-2mm}
 \end{figure}
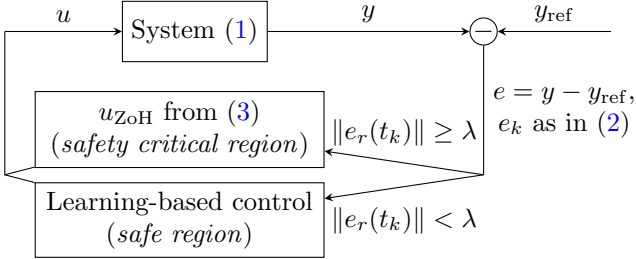

These ideas are based on the observation made in \Cref{Rem:u_eq_0_not_nec}, namely if $\|e_r(t_k)\| < \lambda$, then any bounded input~$u$ can be applied to the system.
In particular, data-driven control schemes are applicable, which often show superior performance due to collection of ``system knowledge'' in terms of input-output data.
The idea of a combined control scheme is illustrated in \Cref{Fig:Controller}. 

Since the calculations in the proof of \Cref{Thm:recursive_feasible_r} involve worst case estimates, the application of $u(t) \neq 0$ for $t \in [t_k, t_k + \tau)$, if $\|e_r(t_k)\| < \lambda$ requires adaption of the sampling time~$\tau$.
This adaption is formulated in the following feasibility result for the switched control strategy
    \begin{equation} \label{eq:combined_controller}
   \! \! \forall  t \in  [t_k, t_{k+1})  :  u(t) =
    \begin{cases}
   \quad  u_{\rm data}, \!\! &  \|  e_r(t_k)\| < \lambda, \\
         - \beta \tfrac{e_r(t_k)}{\|e_r(t_k)\|^2}  , \!\! &\|  e_r(t_k)\| \ge \lambda.
    \end{cases}
\end{equation}
\begin{theorem} \label{Thm:Combined_Controller}
Given a reference {$y_{\rf} \in W^{r,\infty}(\Rp,\R^m)$} and a funnel function $\vp \in \cG$ consider a system~\eqref{eq:Sys_r} with $(d,f,g,\oT) \in \cN^{m,r}$.
Let the constants given in~\eqref{eq:ve_mu_gam}, and $\kappa_0, \kappa_1$ and $\beta$ be given as in \Cref{Thm:recursive_feasible_r}.
Assume that the initial condition satisfies $\chi(y^0)(0) \in \cD_{0}^r$
and, for an activation threshold ${\lambda \in (0,1)}$, and $u_{\rm max} \ge 0$ let the sampling time satisfy
\footnote{
In the published version \url{https://doi.org/10.1016/j.sysconle.2024.105892}, there is a typo: $\lambda^2$ is missing in~\eqref{eq:tau_adapted}.
}
\begin{equation} \label{eq:tau_adapted}
    0 < \tau \le \min \left\{ \frac{ {\lambda^2} \kappa_0}{\kappa_1^2}, \frac{1-\lambda}{\kappa_0 + \|\vp\|_\infty \gM u_{\rm max}} \right\}.
\end{equation}
If $\|u_{\rm data}\|_\infty \le u_{\rm max}$,
then the combined controller~\eqref{eq:combined_controller}
applied to a system~\eqref{eq:Sys_r} yields that $\| e_k(t)\| < 1$ for all $k=1,\ldots,r-1$ and $\|e_r(t)\| \le 1$ for all $t \ge 0$.
This is initial and recursive feasibility of the controller~\eqref{eq:combined_controller}.
In particular, the tracking error satisfies $\|e(t)\| < 1/\vp(t)$ for all $t \in \Rp$.
\end{theorem}
\begin{proof}
    By adapting the sampling time~$\tau$ the statement follows with the same proof as for \Cref{Thm:recursive_feasible_r}.
\end{proof}
With \Cref{Thm:Combined_Controller} at hand, we may now consider the following extensions of the control law~\eqref{eq:controller_recursive}, resulting in a combined controller~\eqref{eq:combined_controller}.

\begin{remark} \label{Rem:AnyControlApplicable}
We comment on some aspects of the presented two-component controller.
\begin{enumerate}[i)]
    \item 
    None of the control schemes applied if $\|e_r(t_k)\| < \lambda$ are required to achieve any tracking guarantees. 
    The only requirement is that the control signal~$u_{\rm data}$ satisfies $\|u_{\rm data}\|_\infty \le u_{\rm max}$ for given $u_{\rm max} > 0$.
    In particular, this means that \emph{any} controller (predictive, or learning-based, or model inversion-based, or locally stabilizing) applied in the safe region satisfies input constraints given by~$u_{\rm max}$.
    Moreover, a control scheme applied in the safe region
    is not even supposed to be suitable for the system to be controlled.
    The latter means that it is possible to apply, e.g., controllers designed for discrete-time systems to the continuous-time system to be controlled.
    Maintenance of the tracking behavior is still ensured by \Cref{Thm:Combined_Controller}.
    \item
     The input $u_{\rm data}$ in~\eqref{eq:combined_controller} is not necessarily supposed to be of data-driven or learning-based type.
    A sample-and-hold version of the funnel control law~\cite{BergIlch21}, i.e.,
    \begin{equation}
        u_{\rm data}(t) = - \alpha(\|e_r(t_k)\|^2) e_r(t_k), \ t \in [t_k,t_{k} + \tau)
    \end{equation}
    is feasible with $u_{\rm max} = \lambda/(1-\lambda^2)$.
    This choice approximates the continuous funnel control signal on a fixed time grid.
    Since this discrete-time funnel controller is safeguarded by the ZoH controller in~\eqref{eq:combined_controller}, none of the issues regarding feasibility of this sampled-and-hold funnel control signal (cf. the considerations in~\cite{BergKaes20}) are present.
    \item
    If a nominal model of the system is available, another combined controller strategy would be to include a pre-computed feedforward signal, cf.~\cite{berger2019combined,drucker2024experimental}, with $u = u_{\rm feedforward} + u_{\rm ZoH}$ where the feedforward controller is active in the safe as well in the safety-critical region. 
    The controller~\eqref{eq:combined_controller} would interpret this additional signal as a ``helpful'' disturbance (``helpful'' since it will reduce the control effort of the feedback), and constraint satisfaction is guaranteed. 
    \end{enumerate}
\end{remark}

\subsection{Data-driven MPC using Willems' fundamental lemma} \label{Sec:WillemsMPC}
In this section, we present a safe region control strategy employing a data-driven MPC scheme. This approach is based on the fundamental lemma by Willems' et al.~\cite{WRMDM05}, which leverages the standard MPC algorithm to a data-enabled predictive control scheme, cf.\ \cite{berberich2022linear,coulson2019data}. Recently, this method attracted a lot of interest in the data-driven control community with various applications e.g.\ in power systems \cite{Schmitz22,huang2019data} and aerial robotics \cite{coulson2019data,elokda2021data}.

Consider a surrogate model for the system~\eqref{eq:Sys_r}. 
The surrogate is given by a discrete-time linear time-invariant system in minimal, i.e.\ controllable and observable, state-space realization
\begin{subequations}
\label{eq:lti_surrogate}
\begin{align}
    x_{k+1} = Ax_k + B u_k\\
    y_k = Cx_k + D u_k
\end{align}
\end{subequations}
with matrices $A\in R^{n\times n}$, $B\in \R^{n\times m}$, $C\in\R^{m\times n}$ and $D\in\R^{m\times m}$. Except the dimension $m$, which is determined by the input and output dimension of the system~\eqref{eq:Sys_r}, the parameters $A,B,C,D$ are assumed to be unknown.

Next we recall the property of persistency of excitation and the fundamental lemma for controllable systems by Willems et al.\ \cite{WRMDM05}, which are pivotal elements in the subsequent discussion. A sequence $u=(u_k)_{k=0}^{N-1}$ with $u_k\in \R^m$, $k=0,\dots, N-1$, is called \emph{persistently exciting of order $L$}, $L\in \N$, if the Hankel matrix
\begin{equation}
    H_L(u) := \begin{bmatrix}
        u_0 & \dots & u_{N-L}\\
        \vdots & \ddots & \vdots\\
        u_{L-1} & \dots & u_{N-1}
    \end{bmatrix}\in \R^{mL\times (N-L+1)}
\end{equation}
has full row rank.
\begin{lemma}[Fundamental lemma]
\label{lem:fl}
    Let $(\hat u,\hat y)=((\hat u_k)_{k=0}^{N-1}, (\hat y_k)_{k=0}^{N-1})$ be an input-output trajectory of length $N$, $N\in\N$, of the system~\eqref{eq:lti_surrogate} such that $\hat u$ is persistently exciting of order $L+n$, where $L\in\N$ and $n$ is the state dimension of system~\eqref{eq:lti_surrogate}. Then $(u,y)=((u_k)_{k=0}^{L-1}, (y_k)_{k=0}^{L-1})$ is an input-output trajectory of length $L$ of system~\eqref{eq:lti_surrogate} if and only if  there is $\nu\in\R^{N-L+1}$ such that
    \begin{equation}
    \label{eq:hankel_colsp}
    \begin{bmatrix}
        u_{[0,L-1]}\\
        y_{[0,L-1]}
    \end{bmatrix} = \begin{bmatrix}
        H_L(\hat u)\\H_L(\hat y)
    \end{bmatrix}\nu.
    \end{equation}
\end{lemma}
The fundamental lemma allows a complete non-parametric, data-driven description of the system's finite-length input-output trajectories based only on measured input-output data.
\begin{remark}
    Note that persistency of excitation order $\tilde L$ implies persistency of excitation of lower order $L$, $L\leq \tilde L$. This fact might be exploited in situations where the state dimension $n$ of a suitable surrogate model~\eqref{eq:lti_surrogate} is unclear but can be estimated, for instance, from physical interpretations of the underlying system~\eqref{eq:Sys_r}. At worst overestimation of $n$ results in an increased data demand for the signal $(\hat u,\hat y)$, while the representation~\eqref{eq:hankel_colsp} is maintained.
\end{remark}

Next we introduce a data-driven MPC scheme leveraged by the fundamental lemma, cf.~\Cref{lem:fl}. To this end let $(\hat u,\hat y)=((\hat u_k)_{k=0}^{N-1}, (\hat y_k)_{k=0}^{N-1})$ be measured input-output data, where $\hat u$ is persistently exciting of order $L+2n$. %
In every discrete time step $t_k$ we aim to solve the optimal control problem
\begin{subequations}
\label{eq:ocp}
\begin{equation}
\label{eq:ocp1}
    \operatorname*{minimize}_{(u,y,\nu,\sigma)} \sum_{i=k+1}^{k+L} \Bigl(\lVert y_i-y_{\text{ref},i}\rVert_Q^2 + \lVert u_i\rVert_R^2\Bigr) + \lambda_\nu\lVert\nu\rVert^2 + \lambda_\sigma\lVert\sigma\rVert^2
\end{equation}
with $(u,y) = ((u_i)_{i=k-n+1}^{k+L},(y_i)_{i=k-n+1}^{k+L})$ subject to
\begin{align}
\label{eq:ocp2}
    \begin{bmatrix}
        u_{[k-n+1,k+L]}\\
        y_{[k-n+1,k+L]} +\sigma
    \end{bmatrix}& = \begin{bmatrix}
        H_{L+n}(\hat u)\\H_{L+n }(\hat y)
    \end{bmatrix}\nu,\\
    \label{eq:ocp3}
    \begin{bmatrix}
        u_{[k-n+1,k]}\\
        y_{[k-n+1,k]}
    \end{bmatrix} &= \begin{bmatrix}
        \tilde u_{[k-n+1,k]}\\
        \tilde y_{[k-n+1,k]}
    \end{bmatrix},\\
    \label{eq:ocp5}
    \begin{bmatrix}
        1 & \dots & 1
    \end{bmatrix}\nu &= 1,\\
    \label{eq:ocp4}
    \lVert u_i\rVert &\leq u_\text{max},\ i=k+1,\dots, k+L
\end{align}
\end{subequations}
on a finite horizon $L>0$, given a past input-output trajectory $(\tilde u,\tilde y)=((\tilde u_i)_{i=k-n}^{k}, (\tilde y_i)_{i=k-n}^{k})$, where $\tilde u_i=u(t_i)$, $\tilde y_i=y(t_i)$ with $u$, $y$ denote the input and output of system~\eqref{eq:Sys_r}, respectively. The weighting matrices $Q,R\in\R^{m\times m}$ in the stage cost in \eqref{eq:ocp1} are assumed to be symmetric and positive-definite. 
As a key difference to standard MPC the state-space model~\eqref{eq:lti_surrogate} is replaced in the optimal control problem~\eqref{eq:ocp} by the equivalent non-parametric description~\eqref{eq:ocp2} based on \Cref{lem:fl}. The constraint~\eqref{eq:ocp3} serves as initial condition which together with the observability of surrogate model~\eqref{eq:lti_surrogate} imposes alignment on the latent state, i.e.\ $x_{[k-n+1,k]}=\tilde x_{[k-n+1,k]}$ for the state sequences $(x_i)_{i=k-n+1}^{k}$ and $(\tilde x_{i})_{i=k-n-1}^{k}$ corresponding to the input-output trajectories $((u_i)_{i=k-n+1}^{k},(y_i)_{i=k-n+1}^{k})$ and $((\tilde u_i)_{i=k-n+1}^{k},(y_i)_{i=k-n+1}^{k})$. 
In order to take into account possible nonlinearities in system~\eqref{eq:Sys_r} not covered by the surrogate~\eqref{eq:lti_surrogate}, we introduce a slack variable $\sigma\in\R^{(L+n)m}$ with weight $\lambda_\sigma>0$ in the cost and the constraint \eqref{eq:ocp5}, cf.\ \cite{Berberich20,berberich2022linear}. Further, the cost functional in~\eqref{eq:ocp1} involves a regularization in terms of $\nu$ with weighting parameter $\lambda_\nu>0$. Further, we impose input constraints in \eqref{eq:ocp4}. The data-driven MPC scheme is summarized in \Cref{alg:1}. 

In practice the observed past trajectory $(\tilde u,\tilde y)$ sampled from the system~\eqref{eq:Sys_r} up to a certain point in time may serve as source for the data $(\hat u,\hat y)$ deployed in the system description~\eqref{eq:ocp3} via Hankel matrices. With this choice more and more data is available with increasing time and, hence, in this way a higher persistency of excitation order can be achieved. As an extension to the above proposed data-driven MPC strategy one may allow for a prediction horizon $L$, which increases over time whenever the updated data is persistently exciting of sufficient order, cf.~\cite[Sec.~5]{schmitz2023safe}.
An additional countermeasure against a divergence of the data-enabled model described by~\eqref{eq:hankel_colsp} and the underlying system due to nonlinearity is to frequently update the data.
\begin{algorithm}
\caption{Data-driven MPC with error guarantees}\label{alg:cap}
\label{alg:1}
\begin{algorithmic}
\State $PE\gets \text{false}$;
\For{$k=0,1,\dots$}
\State get latest sample point $(\tilde u_k,\tilde y_k)$;%
\State calculate $\|e_r(t_k)\|$;
\If{\textbf{not} $PE$}\Comment{learn the dynamics}
    \State update data $(\hat u,\hat y)$, $\hat u_{k} \gets \tilde u_{k}$, $\hat y_{k} \gets \tilde y_{k}$;
    \If{ $\hat u$ is p.e.\ of order $L+n$}
        \State $PE\gets \text{true}$;
        \State store $\mathcal H_{L+n}(\hat u)$, $\mathcal H_{L+n}(\hat y)$;
    \EndIf
\EndIf
\If{$\Norm{e_r(t_k)}<\lambda$}
    \If{PE}
    \Comment{MPC feedback}
        \State $u_\text{act} \gets$ solve(OCP~\eqref{eq:ocp});
    \Else
    \Comment{random input action}
        \State $u_\text{act} \gets \text{random (bounded by $u_{\max}$)}$;
    \EndIf
\Else
\Comment{sampled-data feedback}
    \State $u_\text{act} \gets -\beta \frac{e_r(t_k)}{\Norm{e_r(t_k)}^2}$;
\EndIf
\State apply $u_\text{act}$ as ZoH input action to the system~\eqref{eq:Sys_r}

\EndFor
\end{algorithmic}
\end{algorithm}

We briefly discuss one possible adaption of the previously presented data-driven control algorithm.
The feedback law~\eqref{eq:controller_recursive} involves the recursively defined auxiliary error variables~$e_j$ defined in~\eqref{eq:ek}, which in particular involve higher-order derivatives of both the system output~$y$ and the reference signal~$y_{\rm ref}$. To take the structure of these~$e_j$ into account in the data-driven MPC scheme, we 
discuss one possibility
to include information on these derivatives in the cost function to improve the predictions.
Since the data-driven framework is formulated for discrete-time models~\eqref{eq:lti_surrogate}, we use finite differences to approximate the output's derivatives, i.e., we use $ \frac{y_i - y_{i-1}}{\tau} =: y^{[1]}_i$. Higher-order derivatives are approximated accordingly, and we denote with $y^{[\ell]}_i = \tfrac{1}{\tau^\ell}\sum_{j=0}^\ell (-1)^j \binom{j}{\ell} y_{i-j}$ for~$y_i$ being the output of~\eqref{eq:lti_surrogate} the backwards finite difference approximation of the $\ell^{\rm th}$-order derivative.
Furthermore, we want to take into account the weighting of the higher-order derivatives.
To see, how the derivatives are to be weighted, we explicate the error variable~$e_3$ (we omit the time argument) using the bijection~$\alpha(s) = 1/(1-s)$, and obtain
\begin{equation} \label{eq:e3}
\begin{aligned}
    e_3 &= \vp \ddot e + \frac{1}{1-\|e_2\|^2} e_2 \\
    &= \vp \ddot e + \frac{1}{1-\|e_2\|^2} \Big(\vp \dot e + \frac{1}{1-\|e_1\|^2} e_1 \Big) \\
    &= \vp \Big( \ddot e 
    + \underbrace{\frac{1}{1-\|e_2\|^2}}_{\ge 1}  \dot e 
    + \underbrace{\frac{1}{1-\|e_2\|^2}}_{\ge 1} \underbrace{\frac{1}{1-\|e_1\|^2}}_{\ge 1} e \Big).
\end{aligned}
\end{equation}
From this it is clear that the weighting is decreasing with increasing order of the derivative. 
Combining the regularisation in~\eqref{eq:ocp} and the previous reasoning, we propose the following cost functional
\begin{equation} \label{eq:Willems_Costs_Extended}
\begin{split}
    \sum_{i=k+1}^{k+L} &\Bigl(\sum_{\ell=0}^{r-1} \vp(t_i) \mu_\ell \|y_i^{[\ell]} - y_{\rm ref}^{(\ell)}(t_i)\|^2_Q +  \|u_i\|^2_R\Bigr)\\
    &{ }+ \lambda_\nu\lVert\nu\rVert^2  + \lambda_\sigma\lVert\sigma\rVert^2,
    \end{split}
\end{equation}
where $\mu_0 \ge \mu_1 \ge \cdots \ge \mu_{r-1} \ge 0$, and $\vp(t_i)$ is the funnel function evaluated at~$t=t_i$.
The weights~$\mu_\ell$ reflect the weighting structure in the auxiliary error variables, see~\eqref{eq:e3}.
We observe $1/(1-s^2) = 1$ if and only if~$s=0$, i.e., it is reasonable to order the factors~$\mu_\ell$ strictly.
The reasoning presented above is just one possibility to improve the prediction of the data-enabled MPC
by taking into account the structure of the auxiliary signals~$e_j$ introduced in~\eqref{eq:ek}.
Since the control process is safeguarded by the ZoH controller~\eqref{eq:combined_controller}, there are several options to adapt the cost function in \Cref{alg:1}.

\begin{figure} [htb]
    \centering
    \includegraphics[scale=0.49]{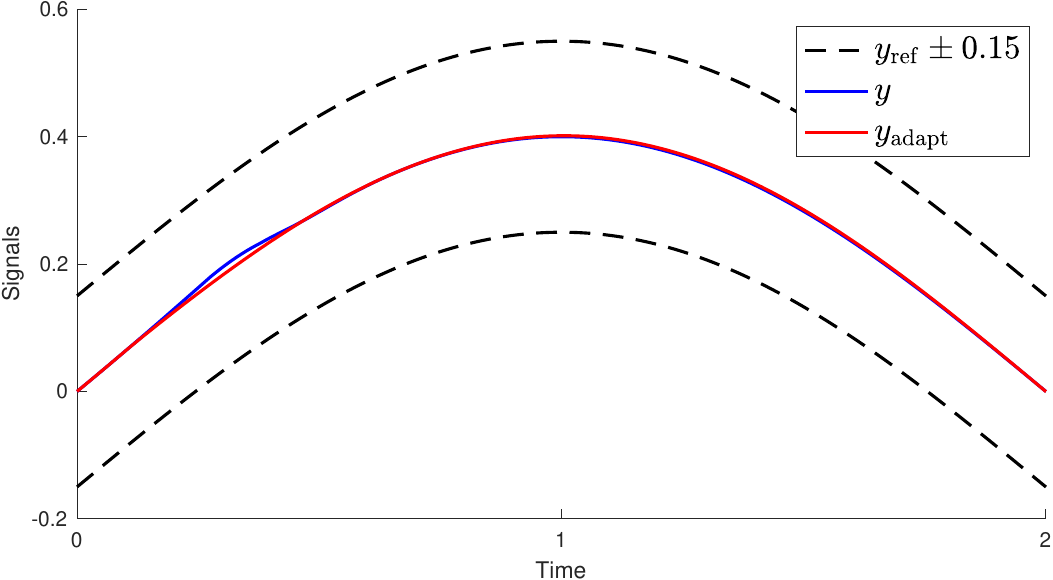}
    \vspace{-2mm}
    \caption{Outputs, reference, and boundaries.}
    \label{fig:dMPC_output}
    \vspace{-5mm}
\end{figure}

\begin{figure} [htb]
    \centering
    \includegraphics[scale=0.49]{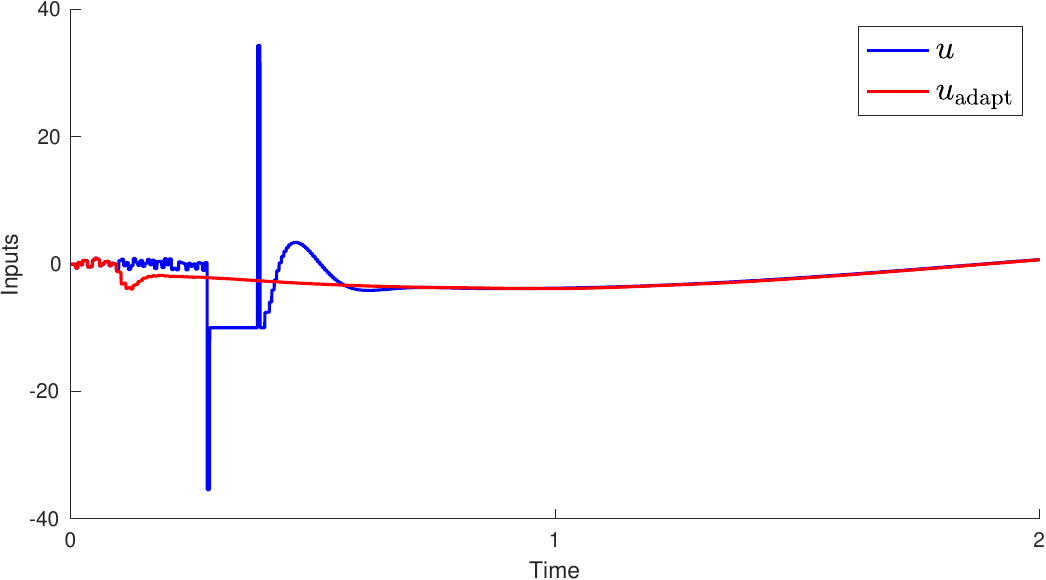}
    \vspace{-2mm}
    \caption{Controls.}
    \label{fig:dMPC_input}
    \vspace{-4mm}
\end{figure}

\begin{figure} [htb]
    \centering
    \includegraphics[scale=0.49]{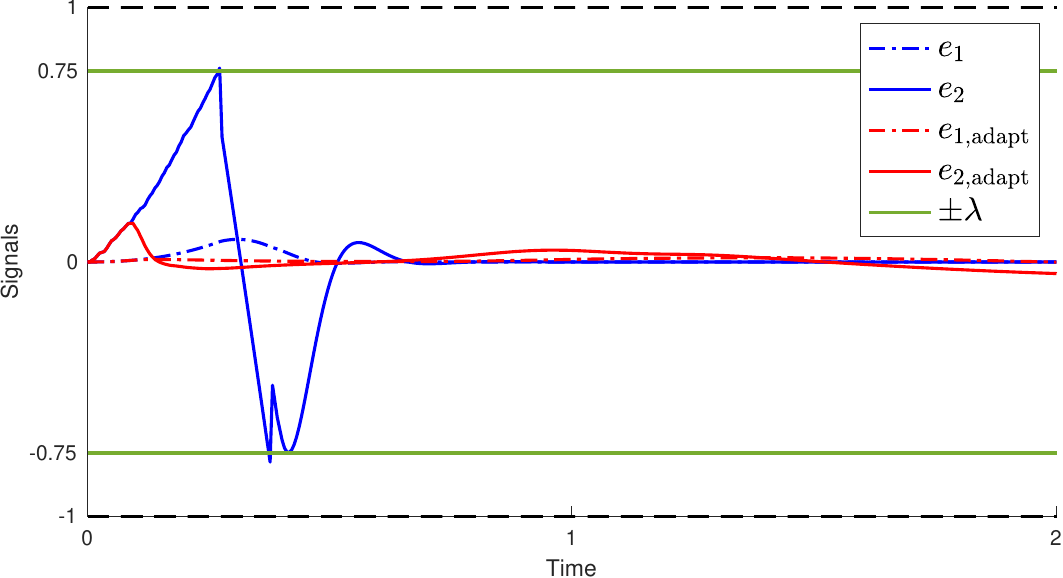}
    \vspace{-5mm}
    \caption{Error variables.}
    \label{fig:dMPC_error}
    \vspace{-4mm}
\end{figure}

In the following we demonstrate the data-enabled MPC scheme described in \Cref{alg:1} on the example system~\eqref{eq:MOC} with fixed prediction horizon $L=20$. Because of the linearity of system~\eqref{eq:MOC} we waive the slack variable 
in the optimal control problem~\eqref{eq:ocp}, i.e. we set $\sigma=0$, and the constraint~\eqref{eq:ocp5}. 
We set $u_\text{max}=10$ which yields~$\tau \le 2.8\cdot 10^{-3}$ according to~\eqref{eq:tau_adapted}. 
As weights we choose $Q=10^{3}\cdot I$, $R=10^{-4}\cdot I$, $\lambda_\nu = 10^{-6}$. We consider a constant funnel given by $\phi(t)=0.15$. The output tracking, the control signal and the auxiliary error variables are depicted in~\Cref{fig:dMPC_output}, \Cref{fig:dMPC_input} and \Cref{fig:dMPC_error} in blue, respectively. In the beginning, there is random control in order to generate a persistently exciting input signal. Then, at $t=0.2728$ persistency of excitation is reached and MPC produces a control signal, however, the error $e_2$ exceeds the safety region~$[- \lambda, \lambda]$. Hence, the ZoH signal becomes active. In the subsequent phase the system is governed by the MPC component, while the signal is saturated at~$-u_\text{max}$. Again the error variable $e_2$ leaves the the safety region at $t=0.3770$ and the ZoH component takes over, resulting in a large control input, which is applied for one
sampling interval. Afterwards, MPC again is sufficient to keep $e_2$ and $e_1$ below $\lambda$ and maintains the tracking goal.

In a second numerical experiment we extend the MPC strategy towards higher auxiliary error variables in the cost functional and an adaptively increasing prediction horizon. 
Respective quantities in \Cref{fig:dMPC_output,fig:dMPC_input,fig:dMPC_error} are labeled with the subscript ``adapt''.
The performance is depicted in \Cref{fig:dMPC_output,fig:dMPC_input,fig:dMPC_error} in red. Starting with $L=1$ the prediction horizon is allowed to increases over time until $L=20$. Further, we set $Q=10^{3}\cdot I$, $R=10^{-4}\cdot I$, $\lambda_\nu = 10^{-6}$ as before, and $\mu_0 = \frac{1}{\phi(0)}$, $\mu_1=\frac{1}{\phi(0)}\cdot 10^{-2}$, where the funnel is constant with $\phi(t)=0.15$. In comparison to the first experiment one observes that the enhanced MPC strategy suffices to safeguard both error variables and, therefore, at no time the ZoH component becomes active. The tracking performance in both runs is of similar quality.

\subsection{Reinforcement Learning: \texorpdfstring{$Q$}{Q}-table control} \label{Sec:RL_Control}
Using the example of $Q$-learning, we show, in this section, how the controller~\eqref{eq:controller_recursive} can be combined 
with model-free Reinforcement Learning~(RL) techniques to safeguard the learning process on the one hand, and to improve the control signal using the control strategy~\eqref{eq:combined_controller}.
$Q$-learning was first developed in~\cite{watkins1989learning}
and has since become a cornerstone of Reinforcement Learning and foundation for many other learning algorithms~\cite{jang2019q}. 

To explain the basic concepts of $Q$-learning,
we consider a nonlinear discrete-time control system of the form
\begin{equation} \label{eq:Sys_RL}
x_{k+1} = f \big(x_k,  u_k \big) 
\end{equation}
where $x\in\cX\subset\R^n$ is the state of the system,
$u\in\cU\subset\R^m$ is the control input,
and $f:\cX\times\cU\to\cX$ is an \emph{unknown} function.
Given an initial state $x^0\in\cX$, we denote, for a control sequence $u=(u_k)\in \cU^\N$,
the solution of~\eqref{eq:Sys_RL} by $x(\cdot;x^0,u)$.
We further assume that there exists a bounded function $r:\cX\times\cU\to\R$,
which is also called~\textit{reward function}.
Note that we do not assume the function~$r$ to be known but merely that 
the reward $r(x_k,u_k)$ can be obtained at every step $k\in\N$ of the system~\eqref{eq:Sys_RL}.
The objective is to maximise the cumulative future reward, i.e. to solve the optimisation problem
\begin{equation}\label{eq:LearningOP}
\mathop {\operatorname{maximize}}_{u\in \cU^\N}  \
             \sum_{k=0}^\infty \gamma^k r(x(k;x^0,u), u_k)
\end{equation}
with discount factor $\gamma\in(0,1)$ which determines the relative importance of long-term versus short-term future rewards.
The so called \textit{$Q$-function} $Q:\cX\times\cU\to\R$, defined by 
\begin{equation} \label{eq:QFunction}
Q(\hat{x},\hat{u}) := r \big(\hat{x},  \hat{u} \big)  
                    +\gamma\sup_{u\in \cU^\N}\sum_{k=0}^\infty \gamma^k r(x(k;f(\hat{x},\hat{u}),u), u_k),
\end{equation}
plays a key role in solving the optimization problem, 
as stated in the following result which we recall for the sake of completeness.
\begin{theorem}[{\cite[Sec. 1.1]{bertsekas2019reinforcement}}]\label{Thm:QPolicyOptimal}
    Consider the system~\eqref{eq:Sys_RL}. 
    If $\pi:\cX\to\cU$ is a feedback control with 
    \begin{equation} \label{eq:OptimalPolicy}
        \pi(x)\in\mathop{\operatorname{argmax}}_{u\in\cU} Q(x,u)
    \end{equation}
    for all $x\in\cX$, then $\pi$ applied to the system~\eqref{eq:Sys_RL} is a solution to the optimization problem~\eqref{eq:LearningOP}.
\end{theorem}

If the Q-function is known, then an optimal feedback control~$\pi$,
in the sense of solving the optimization problem~\eqref{eq:LearningOP}, can be calculated.
Its simplicity makes the  optimal feedback control, also known as the \textit{optimal policy}, very appealing.
This, however, gives rise to the problem of approximating or learning the $Q$-function~\eqref{eq:QFunction}.
While there exist various modern approaches addressing the problem, see e.g.~\cite{jang2019q},
the original $Q$-learning algorithm from~\cite{watkins1989learning} takes the form of~\Cref{alg:Q-learning}.

\begin{algorithm} \caption{Q-learning algorithm} \label{alg:Q-learning}
\begin{enumerate}[1.]
    \item Initialise $k=0$ and $\tilde{Q}_0(x,u):=0$ for all $x\in\cX$ and $u\in\cU$. 
    Let state $x_{0} \in \mathcal{X}$ and select a \textit{learning rate}~$(\alpha_k)\in[0,1)^{\N}$.
    \item Select $u_{k} \in \mathcal{U}$, observe $x_{k+1} = f(x_{k},u_{k}) \in \mathcal{X}$.
    \item Update $\tilde{Q}_{k+1}(x_{k},u_{k})$ by 
    \begin{small}
    \[
    \hspace{-4mm} (1-\alpha_k) \tilde{Q}_k(x_{k},u_{k}) + \alpha_k\!\!
    \left(\!\!r(x_{k},u_{k}) \!+\! \gamma \max_{u' \in \cU}\tilde{Q}_{k}(x_{k+1}, u')\!\!\right)\!\!.
    \]
    \end{small}
    \item Increase $k$ by one, and go to step~2.
\end{enumerate}
\end{algorithm}

An essential part of \Cref{alg:Q-learning} is the selection of the control action in Step~2.
One has to find a balance between selecting the currently expected optimal control and selecting
a different action hoping it yields a higher cumulative reward in the future.
There exist several strategies to address this exploration-exploitation dilemma, see e.g.~\cite{Tijsma16}.
One of the commonly used selection strategies for the control action in the Step~2 of \Cref{alg:Q-learning} is the \textit{$\varepsilon$-greedy} choice.
For a given $\varepsilon\in [0,1]$, the control action is selected as $u_{k+1} = \max_{u \in \cU } \tilde{Q}_{k}(x_{k}, u)$ with probability $1-\varepsilon$, 
and an arbitrary control $u_{k+1}\in\cU$ is selected with probability~$\varepsilon$.

The learning rate~$(\alpha_k)$ also plays a crucial role in addressing the exploration-exploitation dilemma.
It determines the extent to which~\Cref{alg:Q-learning} updates its estimate of the $Q$-function during each iteration by new information.
It is a decisive factor in the convergence rate of the learning algorithm, see e.g.~\cite{even2003learning}.
To proceed combining the controller~\eqref{eq:controller_recursive} with the Q-learning strategy, we recall the following result~\cite{watkins1992q}.

\begin{theorem}[\cite{watkins1992q}]\label{Thm:ResultQLearning} 
    Consider the system \eqref{eq:Sys_RL} with finite sets~$\cX$, $\cU$.
    If the learning rate~$(\alpha_k)\in\ell^2(\N) \backslash \ell^1(\N)$ and if all $(x,u)\in\cX\times\cU$ appear infinitely often in Step~$2$ of~\Cref{alg:Q-learning}, then 
    \[
        \lim_{k\to\infty}\tilde{Q}_k(x,u)=Q(x,u)
    \]
    for all $x\in\cX$, $u\in\cU$.
\end{theorem}
In view of~\Cref{Thm:ResultQLearning}, combining $Q$-learning with the controller~\eqref{eq:controller_recursive}
in the form of a combined controller~\eqref{eq:combined_controller} and applying it to the system~\eqref{eq:Sys_r}
faces three challenges which need to be addressed:
$Q$-learning is formulated for discrete systems,
the sets~$\cX$, $\cU$ are assumed to be finite,
and the problem is presumed to be time-invariant.
Under the assumption that the operator~$\oT$ does not have a time-delay,
using a sampling rate~$\tau>0$ and only applying constant control signals between two sampling instances 
puts the system~\eqref{eq:Sys_r} via evaluation of its solution operator into a discrete system of the form~\eqref{eq:Sys_RL}.
There are various approaches to overcome the requirement of a finite state~$\cX$ and control space~$\cU$, see e.g.~\cite{gaskett1999q}.
As a consequence of~\Cref{Lemma:e_k}, the system states $\chi(y)$, respectively the error signals $e_i$ for $i=1,\ldots, r-1$, 
evolve within a compact set~$K$ when applying the combined controller~\eqref{eq:combined_controller} to the system~\eqref{eq:Sys_r}.
Using a quantization of this compact set therefore is a straightforward way to overcome the problem of the requirement of a finite set~$\cX$.
Since the controller~\eqref{eq:controller_recursive} is bounded by $\beta/\lambda$, a quantization of the set~$\bar{\cB}_{\beta/\lambda}$ is a natural choice for~$\cU$.
However, the curse of dimensionality renders a quantization approach unsuitable for high-dimensional problems.
Note that the quantization of the set $\cU$ is only used for
the learning-based component of the controller~\eqref{eq:combined_controller}
but not for the ZoH controller component~\eqref{eq:controller_recursive} which is used in the safety critical region.
It is still an open research question whether safety guarantees as in~\Cref{Thm:Combined_Controller} 
can be given if the controller~\eqref{eq:controller_recursive} can only emit finitely many different control signals.
Due to the fact that $y_{\rf}$ and~$\phi$ are explicit functions of time, the considered control problem is inherently time variant. 
There are a number of different results for addressing this issue, see e.g.~\cite{hamadanian2022demystifying,khetarpal2022towards}.
Furthermore, it is also possible to encode the time dependency in the state of the system~\eqref{eq:Sys_RL} by enlarging the compact set~$K$ and modifying~\eqref{eq:Sys_RL},
because the functions~$y_{\rf}$ and~$\phi$ are bounded. However, one cannot guarantee that all $(x,u)\in\cX\times\cU$ appear infinitely often in the algorithm,
unless~$y_{\rf}$ and~$\phi$ are periodic. Moreover, encoding the time dependency in the compact set~$K$ further worsens the problem of the curse of dimensionality.
Nevertheless, in virtue of \Cref{Rem:AnyControlApplicable} it is still meaningful to combine the Q-learning scheme with the ZoH controller~\eqref{eq:controller_recursive}.

\begin{figure} [!ht]
    \centering
    \includegraphics[scale=0.55]{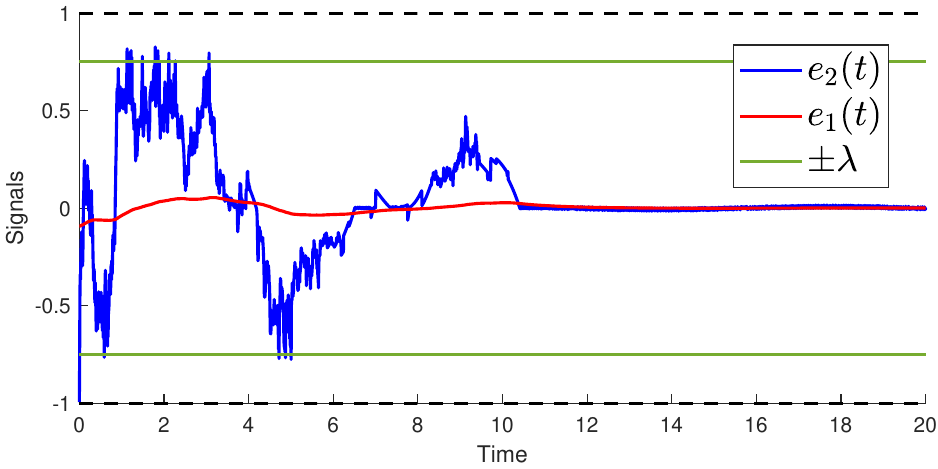}
    \vspace{-2mm}
    \caption{Error signals.}
    \label{fig:Errors_RL}
    \vspace{-5mm}
\end{figure}

\begin{figure} [!ht]
    \centering
    \includegraphics[scale=0.55]{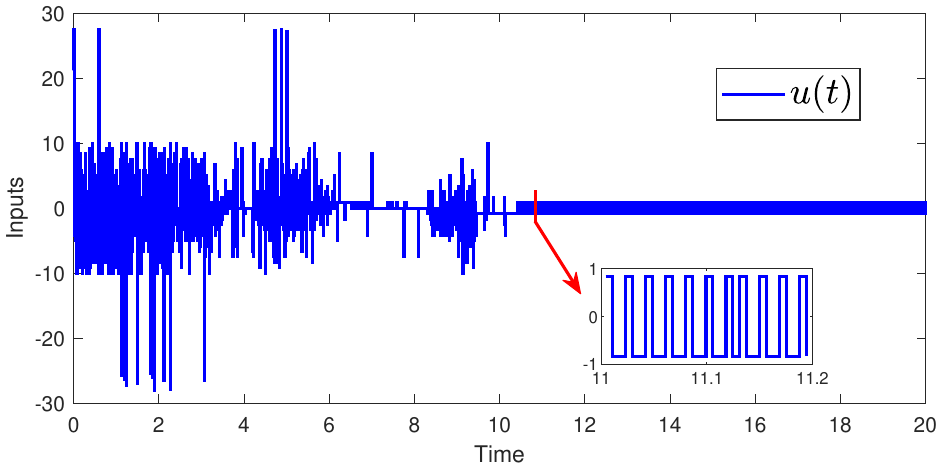}
    \vspace{-2mm}
    \caption{Control signals.}
    \label{fig:Control_RL}
    \vspace{-4mm}
\end{figure}

In the following, we demonstrate the combined controller~\eqref{eq:combined_controller} consisting of~\eqref{eq:controller_recursive} and the $Q$-learning~\Cref{alg:Q-learning} 
on the example system~\eqref{eq:MOC}.
Using the control strategy~\eqref{eq:combined_controller} with sampling time~$\tau>0$ and time instances $t_k\in\tau\N$,
the aim is to take advantage of $Q$-learning by exploring the safe tracking region, e.g. for $\Norm{e_r(t_k)}<\lambda$, 
and applying an improved control signal while the safety critical region is secured by the controller $u_{\rm ZoH}$ as in~\eqref{eq:controller_recursive}
for $\|e_r(t_k)\| \geq \lambda$.
We, therefore, only consider the error variable~$e_r$ for the $Q$-learning \Cref{alg:Q-learning} and
choose a uniform quantization of the set $\bar{\cB}_{\lambda}$ as the state space~$\cX$.
Considering the system~\eqref{eq:Sys_r} and the error variables~\eqref{eq:ek}, $e_r$ satisfies the ordinary differential equation
\begin{align*}
     \dot e_r(t) &= \frac{\dot \vp(t)}{\vp(t)} (e_r(t) - \gamma_{r-1}(t)) + \dot \gamma_{r-1}(t) \\
     &\phantom{=}+ \vp(t) ( f({z}(t))  +g({z}(t)) u - y_{\rf}^{(r)}(t)), 
\end{align*}
with $\gamma_{r-1}(t) := \alpha(\|e_{r-1}(t)\|^2) e_{r-1}(t)$ and ${z}(\cdot) := (d(\cdot), \oT(\chi(y))(\cdot))$.
Sampling this differential equation with sampling time~$\tau$ results in a discrete-time control system.
However, note that it does not have the form~\eqref{eq:Sys_RL} due to the time dependency of $y_{\rf}$ and $\vp$.
Note further that the state variables $e_1,\ldots, e_{r-1}$ are neglected.
Nevertheless, the application of the $Q$-learning algorithm achieves that the error variable~$e_r$ remains, after an initial learning period,
below the threshold $\lambda$ as simulations show, see~\Cref{fig:Errors_RL,fig:Control_RL}.
Further research is necessary to determine whether it is always the case that solely considering $e_r$ in the $Q$-learning algorithm is sufficient
and if guarantees about the convergence of the learning algorithm can be given despite the inherent time dependency of the problem.
As for the set of control values, we choose $\cU$ to be a uniform quantization of the set~$\bar{\cB}_{u_{\max}}$ where $u_{\max}$ is chosen as $u_{\max}=10$ as in the example in \Cref{Sec:WillemsMPC}.
To improve the performance of the original controller~\eqref{eq:controller_recursive},
meaning better tracking performance and reduced control values, we choose the reward function 
\[
    r(e_r(t_k),u) = -\Norm{e_r(t_k)}^2 - \alpha_u \Norm{u}^2,
\]
with parameter $\alpha_u\in\Rp$.
The function $r$ rewards small values of the error variable~$e_r$ and the applied control values (depending on the penalty parameter $\alpha_u$).
For the simulation of the example system~\eqref{eq:MOC}, we chose the system parameters as in~\Cref{Sec:Example}. 
The reference trajectory was $y_{\rf} = 0.4 \sin(\tfrac{\pi}{4}t)$ for $t \in [0,20]$.
Further for the $Q$-learning parameters, the size of the finite sets $\cX$ and $\cU$ were selected as $8$ and $25$, respectively.
The learning rate was set as constant $\alpha = 0.8$.
In order to let the algorithm explore the state and action space, the greedy parameter was set to $\varepsilon = 1$ for $t \in [0,1]$,
thereafter the greedy parameter is halved in order to take the control action more often according to learned $Q$-function.
For the reward function the parameter $\alpha_u=1/u_\text{max}$ was selected.
The simulations are depicted in~\Cref{fig:Errors_RL,fig:Control_RL}.
\Cref{fig:Errors_RL} shows how the error signals evolving within the funnel, respectively the $\lambda$ activation threshold.
\Cref{fig:Control_RL} shows the corresponding control action.
It can be seen that with the help of the primary controller $u_{\rm ZoH}$ in~\eqref{eq:controller_recursive},
Q-learning algorithm is able to safely explore the state and action space
and learn/approximate the Q-function by applying random control actions with an amplitude lower than $10$.
Only if the error exceeds the activation threshold,
the ZoH control component intervenes with a large control input to prevent a violation of the funnel boundaries. 
One can see  that with decaying $\varepsilon$ the number of random control actions applied to the system reduces 
and the auxiliary signal $e_2(t)$ gets closer to $0$ and remains close to it.
Overall, the $Q$-learning algorithm reduces the peaks of the control significantly 
in comparison to~\Cref{Sec:Example} where merely the controller~\eqref{eq:controller_recursive} was applied.

\begin{remark}
    To reduce computational effort, the control signal $u_{\rm data}$ in~\eqref{eq:combined_controller} does not have to be updated at every~$t_i=i \tau$. 
    Since the system class~\eqref{eq:Sys_r} allows for bounded disturbances, it is possible to combine the data-driven control with a move blocking strategy, cf.~\cite{cagienard2007move}, i.e., to apply the control value~$u_{\rm data}$ for longer than one sampling interval~$\tau$.
    If then $e_r$ leaves the safe region, the controller~\eqref{eq:combined_controller} interprets the additional value~$u_{\rm data}$ as a disturbance in the system (according to \Cref{Ass:d_bounded} this means $D=\|d\|_\infty + u_{\rm max}$), and hence the constraint satisfaction is guaranteed by the controller.
    Note that system measurements, however, have to be taken at every~$t_i = i \tau$.
\end{remark}

\subsection{Nonlinear example: Van der Pol oscillator}
\begin{figure} [!ht]
    \centering
    \includegraphics[scale=0.47]{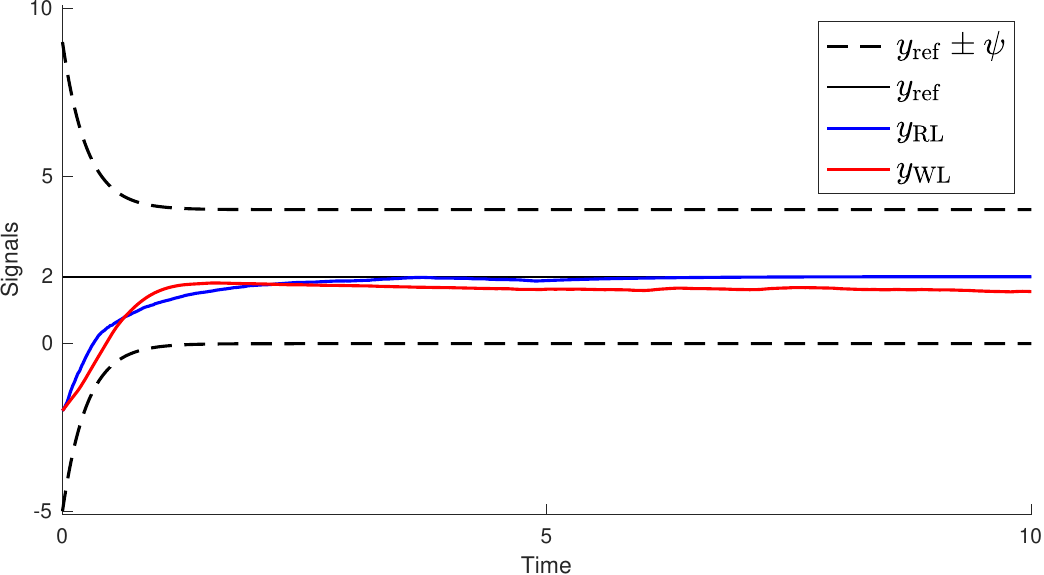}
    \vspace{-2mm}
    \caption{Output within the funnel around reference.}
    \label{fig:Errors_VanDerPol}
    \vspace{-5mm}
\end{figure}

\begin{figure} [!ht]
    \centering
    \includegraphics[scale=0.47]{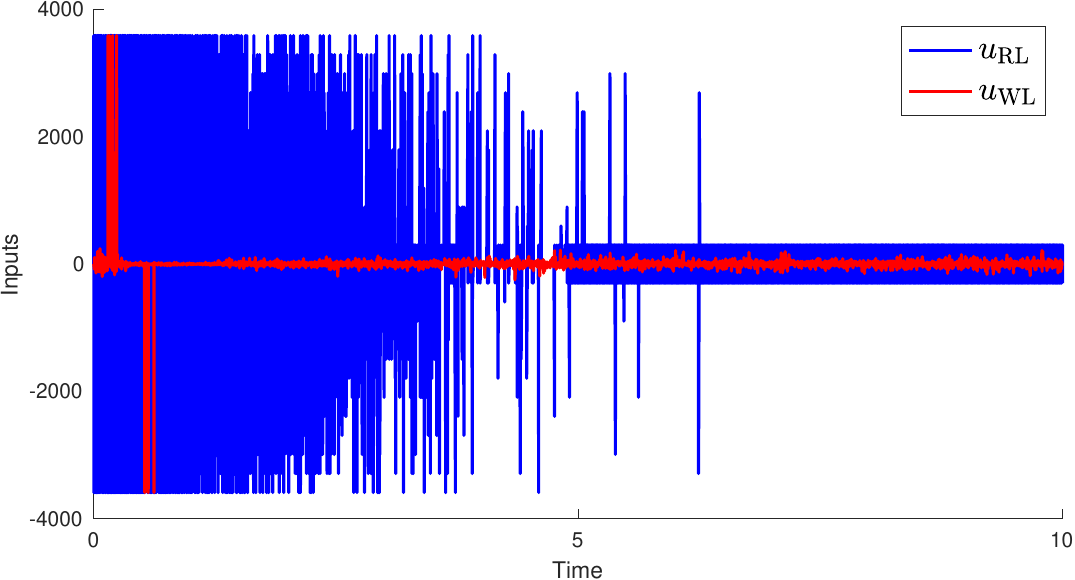}
    \vspace{-2mm}
    \caption{Control signals.}
    \label{fig:Control_VanDerPol}
    \vspace{-4mm}
\end{figure}

The previously considered example~\eqref{eq:MOC} is a linear system.
Now, we briefly present a numerical simulation of a nonlinear system, namely an externally driven Van der Pol oscillator with additive disturbance~$d(t)$.
This system is a typical example for nonlinear systems with global relative degree two.
The system dynamics are governed by
\begin{equation*}
    \ddot y(t) - (1- y(t)^2) \dot y(t) + y(t) - u(t) - d(t) = 0,  
\end{equation*}
with $y(0), \dot y(0) \in \R$, and external input~$u(t)$. 
The bounded function~$d \in L^\infty(\Rp,\R)$ acts as a disturbance.
To illustrate the effect of the nonlinearity, we track a constant reference~$y_{\rm ref}(t) \equiv 2$. 
We use a non-constant funnel $\vp(t) = (a e^{-bt} + c)^{-1}$ with $a = 5, b = 4, c=2$.
We set $y(0) = -2$, i.e., starting with a large tracking error, and $\dot y(0) = 4$.
We choose $d(t) = 0.1 \cos(7t)$.
For activation threshold $\lambda = 0.75$ we calculate $\beta \ge 2.6918 \cdot 10^3$.
Then, with $u_{\rm max} = \beta/\lambda$ (input bound for data-driven control) the requirements of \Cref{Thm:Combined_Controller} are satisfied with 
$\tau \le 1.149 \cdot 10^{-4} $.
The results of the simulation are depicted in \Cref{fig:Errors_VanDerPol,fig:Control_VanDerPol},
where we use the subscript WL for using \Cref{alg:1}, and RL for using \Cref{alg:Q-learning}.
Like for the previous example, we consider both data-driven controllers.
For \Cref{alg:1} we use the cost function~\eqref{eq:Willems_Costs_Extended} with parameters
$\phi(t_i)\mu_0(t_i)=1$, $\phi(t_i)\mu_1(t_i) = 1/2 \cdot 10^{-3}$, $Q=2\cdot10^3\cdot I$, $R = 10^{-4}\cdot I$, $\lambda_\nu = 10^{-5}$, $ \lambda_\sigma = 10^6$.
Moreover, we use an adapted prediction horizon with $L_\mathrm{max} = 20$. To account for the underlying nonlinear dynamics, additionally, once $L_\mathrm{max}$ is reached, the database is constantly updated with new data points and older data is removed, provided that the persistency of excitation is maintained.
As can be seen in \Cref{fig:Errors_VanDerPol} the tracking guarantees are valid as stated in \Cref{Thm:Combined_Controller}.
\Cref{fig:Control_VanDerPol} shows the control actions.
It can be observed that the Q-learning algorithm requires some more data to produce controls,
which are sufficient to achieve $\|e_2(t_k)\| < \lambda$, i.e., to avoid activation of the safeguarding controller component.
While \Cref{alg:1} (WL) produces smaller input values, the tracking is more accurate using \Cref{alg:Q-learning} (RL) after the learning process.

\section{Conclusion and future work} \label{Sec:Conclusion}
We presented a novel two-component controller for continuous-time nonlinear control systems.
The ZoH tracking controller consists of a data-driven/learning-based component and a 
discrete-time output-feedback controller with prescribed performance.
The feedback 
controller is designed to achieve the control objective (tracking  
with prescribed performance) 
and safeguards the 
learning-based controller.
We derived explicit upper bounds on the sampling time~$\tau > 0$ and for the maximal control input.
As data-driven controller we employed an MPC algorithm based on the fundamental results of Willems et al.~\cite{WRMDM05}, which enables predictive control using only input-output data.
Further, we implemented a Reinforcement Learning scheme and investigated a Q-table control algorithm to explore the system's dynamics.
The proposed two-component data-driven controller was
proven to achieve the control objective, and in particular, outperform the pure feedback controller.

Based on the presented results, future work will aim to reduce the conservatism of the controller and to investigate the interplay with observers and/or the funnel pre-compensator \cite{BergReis18b,Lanz22} to alleviate the strict assumption of not only {knowing} the output but also its derivatives.
Moreover, we plan to perform a comprehensive comparison (simulation study) with other data-driven 
ZoH controllers, e.g., the one recently proposed in~\cite{BoldGrun23}, and combining these with the proposed safeguarding feedback component.

\bibliographystyle{IEEEtran}
\bibliography{references.bib}

\appendix
\begin{appendices}
  \crefalias{section}{appsec}
\section{Proofs of auxiliary results} \label{Sec:Appendix}

We present the proofs of the auxiliary results \Cref{Lemma:e_k,Lemma:DynamicsBounded} presented in \Cref{Sec:Aux_results}, and \Cref{Thm:recursive_feasible_r} in \Cref{Sec:MainResult}.

\begin{proof}[Proof of \Cref{Lemma:e_k}] 
We use the constants~$\ve_k,\mu_k>0$ defined in~\eqref{eq:ve_mu_gam}, 
and to improve legibility, we use the notation $e_k(t):=e_k(t,\chi(\zeta)(t))$ for $\zeta \in \cY_\delta^r$.
Let $\delta\in(0,\infty]$ and  $\zeta\in \cY^r_\delta$ be arbitrary but fixed.
We define the auxiliary function $\gamma_k(t) := \alpha(\|e_k(t)\|^2) e_k(t)$, and set $\gamma_0(\cdot) = \dot \gamma_0(\cdot) = 0$.
Note that for $k=1,\ldots,r-1$ each of the error signals defined in~\eqref{eq:ek} satisfies for $t\in[0,\delta)$ the differential equation
\[
    \dot e_k 
    = \frac{\dot \vp}{\vp} ( e_k \!-\!\gamma_{k-1} )\! +\! e_{k+1}\! 
    \!+\! \dot \gamma_{k-1}
    - \!\alpha(\|e_k\|^2) e_k ,
\]
where the dependency on $t$ has been omitted and $e^{(k)}$ denotes the $k$-th derivative of $e(t) = \zeta(t)-y_{\rf}(t)$.
We observe 
\begin{align*}
    \dot \gamma_k &= 2 \alpha'(\| e_k\|^2) \al e_k, \dot e_k \ar e_k + \alpha(\| e_k\|^2) \dot e_k.
\end{align*}
Seeking a contradiction, we assume that for at least one ${\ell \in \{1,\ldots,r-1\}}$ there exists $t^* \in (0,\delta)$
such that $\|e_\ell(t^*)\|^2 > \ve_\ell$. W.l.o.g. we assume that this is the smallest possible $\ell$.
Invoking $\chi(y^0)\in\cD_0^r$ and continuity of the involved functions we may define $t_* := \max \setdef{ t \in [0,t^*) }{ \| e_\ell(t) \|^2 = \ve_\ell}$.
Then, for $t \in [t_*,t^*]$ we calculate, omitting again the dependency on~$t$,
\begin{small}
\begin{align*}
    &\dd{t} \tfrac{1}{2} \| e_\ell\|^2 
     =\! \al e_\ell, \tfrac{\dot \vp}{\vp} ( e_\ell - \gamma_{\ell-1} ) + e_{\ell+1}  + \dot \gamma_{\ell-1} - \alpha(\|e_\ell\|^2) e_\ell\ar \\
    & \le\! \| e_\ell \| \!\left( \SNorm{\frac{\dot \vp}{\vp}}\!\!\!\! ( 1 \!+\! \alpha(\ve_{\ell-1}^2) \ve_{\ell-1}) 
    \! + \! 1  \! + \! \bar \gamma_{\ell-1}\!  - \!\alpha(\ve_\ell^2) \ve_\ell\! \right) 
    \le \!0,
\end{align*}
\end{small}
in the last line we used the monotonicity of $\alpha(\cdot)$,
the definition of $\ve_\ell$, and that $\dot \gamma_{\ell -1}$ is bounded by minimality of~$\ell$. 
Hence, the contradiction {$\ve_\ell < \| e_\ell(t^*)\|^2 \leq \| e_\ell(t_*)\|^2 = \ve_\ell$} arises after integration.
This yields boundedness of $e_\ell, \gamma_\ell$. 
Using the derived bounds we estimate
\[
    \| \dot e_\ell \|  \le  \SNorm{\frac{\dot \vp}{\vp} }\!\!\!\! ( 1 + \alpha(\ve_{\ell-1}^2) \ve_{\ell-1} ) 
    + 1  + \alpha(\ve_\ell^2) \ve_\ell 
    + \bar \gamma_{\ell-1}
    = \mu_\ell. 
\]
We conclude
$\| e_k(t) \| \le \ve_k < 1$ and $\| \dot e_k(t)\| \le \mu_k$ for all $k = 1,\ldots,r-2$ and all $t \in [0,\delta)$. 
For $k=r-1$ the same arguments are valid invoking $e_r : [0,\delta) \to \overline{\cB_1}$.
\end{proof}

\begin{proof}[Proof of \Cref{Lemma:DynamicsBounded}]
    To prove the assertion, we invoke continuity of the system functions~$f,g$ and the resulting boundedness on compact sets.
    According to~\Cref{Lemma:e_k}, there exist ${\varepsilon_k \in (0,1)}$ for $k=1,\ldots,r-1$
    such that 
    \[
        \!\fa\! \zeta\!\in\!\cY_{\infty}^r \fa t\in\Rp \!\fa k=1,\ldots, r-1: \Norm{e_k(t,{\chi}(\zeta)(t))}\!{\leq}\varepsilon_k.
    \]
    {Further, $\|e_{r}(t,\chi(\zeta)(t))\|\leq 1$.}
    Thus, due to the definition of~$e_k$ in~\eqref{eq:ek}, there exists a compact set $K_\zeta\subset\R^{rm}$ with
    \[
        \fa\zeta\in\cY_{\infty}^r\fa t\in\Rp:\quad \chi(\zeta)(t)\in K_\zeta.
    \]
    Due to the BIBO property of the operator~$\oT$, there exists a compact set $K_q\subset\R^q$ 
    with $\oT(\xi)(\Rp)\subset K_q$ for all ${\xi\in\con(\Rp,\R^{rm})}$ with $\xi(\Rp)\subset K_{\zeta}$.
    For arbitrary ${\delta\in(0,\infty)}$ and $\zeta\in\cY_{\delta}^r$, we have{, according to~\Cref{Lemma:e_k},}
    \[
        \fa t\in[0,\delta) \fa k=1,\ldots, r-1:\quad \Norm{e_k(t,{\chi}(\zeta)(t))}{\leq}\varepsilon_k.
    \]
    {Further, $\|e_{r}(t,\chi(\zeta)(t))\|\leq 1$.}
    Thus, $\chi(\zeta)(t)\in K_{\zeta}$ for all ${t\in [0,\delta)}$.
    For every element $\zeta\in\cY_{\delta}^r$ the function $\chi(\zeta)|_{[0,\delta)}$
    can smoothly be extended to a function $\tilde{\zeta}\in(\con(\Rp,\R^m))^r$ 
    with $\tilde{\zeta}(t)\in K_{\zeta}$ for all $t\in\Rp$.
    Due to the BIBO property of the operator $\oT$, we have $\oT(\tilde{\zeta})(t)\in K_q$ for all $t\in\Rp$.
    Since~$\oT$ is causal, this implies $\oT(\chi(\zeta))|_{[0,\delta)}(t)\in K_q$ for all $t\in[0,\delta)$ and $\zeta\in\cY_{\delta}^r$.
    Define the compact set $K:=\overline{\cB_D}\times K_{q}\subset\R^{p + q}$.
    Since $f(\cdot)$ and $g(\cdot)$ are continuous,
    the constants $\fM:=\max_{x \in K} f(x)$ and $\gM:=\max_{x\in K} g(x)$ exist.
    For every $\delta\in(0,\infty]$, $\zeta\in \cY_{\delta}^r$, and $d\in L^\infty(\Rp,\R^p)$ with $\SNorm{d}\leq D$ we have
        $\fa t\in[0,\delta):\quad  (d(t),\oT(\chi(\zeta))(t))\in K$.
    Therefore, we obtain $\fM \geq \SNorm{f((d,\oT(\chi(\zeta)))|_{[0,\delta)})}$ and ${\gM \geq \SNorm{g((d,\oT(\chi(\zeta)))|_{[0,\delta)})}}$.
    Since $g(x)$ is positive definite, for every $x\in K$ there exists $\gm>0$ such that 
    $ \gm \leq \frac{\al z, g((d,\oT(\chi(\zeta)))|_{[0,\delta)}(t))z\ar }{\Norm{z}^2}$ for all $z\in\R^m\backslash\cbl0\cbr$.
\end{proof}

\begin{proof}[Proof of \Cref{Thm:recursive_feasible_r}] 
The proof consists of two main steps. 
In the first step we establish the existence of a solution of the initial value problem~\eqref{eq:Sys_r},~\eqref{eq:controller_recursive}.
In the second step we show feasibility of the proposed control law, i.e., all error variables are bounded by $\ve_k$ and the tracking error evolves within the funnel boundaries.

\noindent
\emph{Step 1.}
The application of the control signal~\eqref{eq:controller_recursive} to system~\eqref{eq:Sys_r} leads to an initial value problem. 
If this problem is considered on the interval~$[0,\tau]$, then there exists a unique maximal solution on $[0,\omega)$ with $\omega\in(0,\tau]$.
If all error variables $e_k$ evolve within the set $\cB_1$ for all $t\in [0,\omega)$, then $\| \chi(y)(\cdot) \|$ is bounded on the interval $[0,\omega)$ and, as a consequence of the BIBO condition of the operator, $\oT(\cdot)$ is bounded as well. 
Then $\omega =\tau$, cf. \cite[\S~10, Thm.~XX]{Walt98} and there is nothing else to show.
Seeking a contradiction, we assume the existence of $t \in [0,\omega)$ such that $\| e_k(t)\| \ge 1$ for at least one $k = 1,\ldots,r$.
Invoking \Cref{Lemma:e_k} it remains only to show that the last error variable $e_r$ satisfies $\| e_r(t) \| \le 1$ for all $t \in [0,\omega)$.
Before we do so, we record the following observation.
For $\gamma_{r-1}(t) := \alpha(\|e_{r-1}(t)\|^2) e_{r-1}(t)$ we calculate for ${z}(\cdot) := (d(\cdot), \oT(\chi(y))(\cdot))$
\begin{equation} \label{eq:J}
\begin{aligned}
    & \dot e_r(t) - \vp(t) g({z}(t)) u = \dot \vp(t) e^{(r-1)}(t) + \vp(t) e^{(r)}(t) \\
    & \quad + \dot \gamma_{r-1}(t) - \vp(t) g({z}(t)) u \\
    &= \frac{\dot \vp(t)}{\vp(t)} (e_r(t) - \gamma_{r-1}(t)) + \dot \gamma_{r-1}(t) \\
     &\quad + \vp(t) ( f({z}(t)) - y_{\rf}^{(r)}(t)  ) =: J(t).
\end{aligned}
\end{equation}
\noindent
\emph{Step 2.} We show $\|e_r(t)\| \le 1$ for all $t \in [0,\omega)$.
We separately investigate the two cases $\| e_r(0)\| < \lambda$ and $\| e_r(0)\| \ge \lambda$. \\
\emph{Step 2.a} We consider $\| e_r(0)\| < \lambda$. In this case we have $u = 0$.
Seeking a contradiction, we suppose that there exists $t^* := \inf \setdef{ t \in (0,\omega)}{ \| e_r(t) \|> 1} $. For the function $J(\cdot)$ introduced in~\eqref{eq:J} we observe $\| J|_{[0,t^*)} \|_\infty \le \kappa_0$ according to \Cref{Lemma:e_k,Lemma:DynamicsBounded}.
Then we calculate for $t \in [0,t^*]$
\begin{align*}
  1 &=  \| e_r(t^*) \| \le \|e_r(0)\|  + \textstyle \int_0^{t^*}\!\!  \| \dot e_r(s) \|\!  \d s \\
    &= \|e_r(0)\| + \textstyle \int_0^{t^*}   \| J(s) \| \d s \\
    &\le  \|e_r(0)\| +  \textstyle \int_0^{t^*} \kappa_0  \d s 
    < \lambda + \kappa_0 \omega < 1,
\end{align*}
where we used $t^* < \omega \le \tau < (1-\lambda)/\kappa_0$.
This contradicts the definition of $t^*$. \\
\emph{Step 2.b}
We consider $\| e_r(0)\| \ge \lambda$. In this case we have the control $u = - \beta e_r(0)/\|e_r(0)\|^2$.
We show again $\|e_r(t)\| \le 1$ for all $t \in [0,\omega)$.
To this end, seeking a contradiction, we suppose the existence of $t^* = \inf \setdef{(0,\omega)}{ \| e_r(t) \| > 1 }$.
    Invoking the initial conditions and continuity of the involved functions, and utilizing~\Cref{Lemma:DynamicsBounded} {and \eqref{eq:J}}, we calculate for $t \in [0,t^*]$
    \begin{align*}
       & \dd{t} \tfrac{1}{2} \| e_r(t)\|^2 = \al e_r(t), \dot e_r(t) \ar 
         = \al\! e_r(0) + \!\! \textstyle \int_0^t\!\! \dot e_r(s) \d s, \dot e_r(t) \!\ar \\
        & \le \|e_r(0)\| \| J(t) \| + \omega \| \dot e_r|_{[0,t^*]}\|^2_\infty + \vp(t) \al e_r(0), g({z}(t)) u \ar \\
        & = \|e_r(0)\| \| J(t) \| \!+\! \omega \| \dot e_r|_{[0,t^*]}\|^2_\infty \!-\! \vp(t) \beta \tfrac{\al e_r(0), g({z}(t)) e_r(0)\ar}{\|e_r(0)\|^2} \\
        & \le \|e_r(0)\| \kappa_0 + \omega \| \dot e_r|_{[0,t^*]}\|^2_\infty - \inf_{s \ge 0} \vp(s) \gm \beta \\
        & \le \kappa_0 + \omega {\kappa_1^2/\lambda^2} - \inf_{s \ge 0} \vp(s) \gm \beta  
         \le 2 \kappa_0 - \inf_{s \ge 0} \vp(s) \gm \beta  < 0,
    \end{align*}
    the third line due to $t^* < \omega \le \tau$, the penultimate line via the definition of $\tau$ and the last line by definition of~$\beta$; moreover, we used $\| \dot e_r|_{[0,t^*]} \| \le \kappa_1$ and $\| J|_{[0,t^*]} \|_\infty \le \kappa_0$.
    In particular this yields $\tfrac{1}{2} \dd{t} \| e_r(t)|_{t=0}\|^2 < 0$, by which $t^* > 0$.
    Therefore, we find the contradiction $1 = \| e_r(t^*)\|^2 < \| e_r(0)\|^2 \le 1$.
    Repeated application of the arguments in \emph{Steps 1} and \emph{2} on the interval $[t_i,t_i + \tau]$, $i \in \N$, yields recursive feasibility.
\end{proof}
\end{appendices}

\end{document}